\pdfoutput=1
\documentclass[12pt]{amsart}%
\usepackage{amscd, amssymb, amsfonts, amsthm, array}
\usepackage{amsmath}
\usepackage{amsfonts}
\usepackage{amssymb}
\usepackage{graphicx}%
\evensidemargin=\oddsidemargin 
\email{glickenstein@math.arizona.edu}
\thanks{DG partially supported by NSF grant DMS 0748283.}
\email{payntrac@isu.edu, tpayne@member.ams.org}
\thanks{TLP partially supported by a WeLead grant from NSF Advance grant RRES31.}
\theoremstyle{plain}
\newtheorem{thm}{Theorem}
\newtheorem{theorem}[thm]{Theorem}
\newtheorem{prop}[thm]{Proposition}
\newtheorem{proposition}[thm]{Proposition}

\newtheorem{lemma}[thm]{Lemma}

\newtheorem*{theorema}{Theorem A}
\newtheorem*{theoremb}{Theorem B}
\theoremstyle{definition}   
\newtheorem{definition}[thm]{Definition}

\theoremstyle{remark}
\newtheorem{remark}[thm]{Remark}
\numberwithin{equation}{section}
\numberwithin{thm}{section}

\begin{document}
\author{David Glickenstein}
\address{Department of Mathematics, University of Arizona, Tucson, AZ 85721}
\author{Tracy L.~Payne}
\address{Department of Mathematics, Idaho State University, Pocatello, ID 83209-8085}
\date{\today}
\title[RF on 3D metric Lie algebras]{Ricci flow on three-dimensional, unimodular metric Lie algebras}

\begin{abstract}
{ }We give a global picture of the Ricci flow on the space of
three-dimensional, unimodular, nonabelian metric Lie algebras considered up to
isometry and scaling. The Ricci flow is viewed as a two-dimensional dynamical
system for the evolution of structure constants of the metric Lie algebra with
respect to an evolving orthonormal frame. This system is amenable to direct
phase plane analysis, and we find that the fixed points and special
trajectories in the phase plane correspond to special metric Lie algebras,
including Ricci solitons and special Riemannian submersions. These results are
one way to unify the study of Ricci flow on left invariant metrics on
three-dimensional, simply-connected, unimodular Lie groups, which had
previously been studied by a case-by-case analysis of the different Bianchi
classes. In an appendix, we prove a characterization of the space of
three-dimensional, unimodular, nonabelian metric Lie algebras modulo isometry
and scaling.

\end{abstract}
\maketitle

%-------------------------------------------------------------------------
%-------------------------------------------------------------------

\section{Introduction}

In this paper we study the Ricci flow on three-dimensional, unimodular metric
Lie algebras. Metric Lie algebras are in one-to-one correspondence with
left-invariant Riemannian metrics on simply-connected Lie groups, and Ricci
flow on such metrics has been studied by a number of authors (e.g., \cite{IJ}
\cite{KM} \cite{IJL} \cite{Lot1} \cite{G2} \cite{P2} \cite{CS2} \cite{Lau4}).
The major advances in this paper are (1) a unification of the trajectories for
the Ricci flow, previously viewed individually in case-by-case studies of
Bianchi classes, into a single global topological picture, and (2) use of a
new technique of flowing the Lie structure constants, which highlights
different features of the system than the usual evolution of metric coefficients.

The space of metric Lie algebras has been studied by a number of authors
(e.g., \cite{Jen} \cite{Lau2} \cite{Lau4}). Understanding Ricci flow on the
space of metric Lie algebras is important for studying both homogeneous spaces
and Ricci flow of general manifolds. A number of Ricci soliton metrics (fixed
points of the Ricci flow up to diffeomorphism invariance and rescaling) have
been found on homogeneous spaces (see, e.g., \cite{BD} \cite{Lau1} \cite{Lau3}
\cite{P1} \cite{Guz} \cite{Lau4}), and it has been suggested that finding
Ricci solitons may be a promising way to attack Alekseevskii's conjecture (see
\cite{Lau3} and \cite{Lau4}). Lott has shown that three-dimensional, Type III
solutions to Ricci flow converge to the known homogeneous expanding solitons
as they collapse in the limit (see \cite{Lot2}). Ricci flow on homogeneous
spaces is also useful in constructing self-dual solutions of Euclidean vacuum
Einstein's equations (see \cite{BEPS}).

We will consider the set $\mathcal{M}$ of three-dimensional, nonabelian,
unimodular metric Lie algebras modulo isometry and scaling. Milnor gives an
excellent description of such metric Lie algebras in \cite{Mil}, in particular
showing that there exists a special orthonormal basis $\left\{  e_{1}%
,e_{2},e_{3}\right\}  $ which diagonalizes both the Ricci endomorphism and the
Lie bracket (we say that the Lie bracket is diagonalized if $\left[
e_{i},e_{j}\right]  $ is a scalar multiple of $e_{i}\times e_{j}$). Thus the
set of three-dimensional, unimodular metric Lie algebras depends only on three
parameters. In fact, there are two natural choices of those three parameters,
and the Ricci flow through these parameter spaces takes one of the following forms:

\begin{enumerate}
\item Fix the Lie algebra and let the metric vary.

\item Evolve the frame to keep it orthonormal and let the structure constants vary.
\end{enumerate}

\noindent In both cases, the Lie bracket and inner product remain diagonal
with respect to the frame. However, in the first case the Lie bracket
coefficients are fixed and the lengths of basis elements change. In the second
case, the Lie bracket coefficients change but the lengths of basis elements do
not change (since the basis evolves to stay orthonormal). It is extremely
important that the frame remains orthogonal under the flow, which follows from
the fact that both the structure constants and the Ricci curvature can always
be diagonalized at the same time as the metric. This is true for
three-dimensional, unimodular metric Lie algebras, but not in general. The
lack of such a frame is the major obstacle for classifying Ricci flow on
four-dimensional, simply-connected homogeneous spaces; Isenberg-Jackson-Lu
\cite{IJL} classify Ricci flow for some Riemannian homogeneous spaces which do
admit such a frame.

Since $\mathcal{M}$ is a three-dimensional space considered modulo rescaling,
we have a two-dimensional system of ODEs, which is reasonable to analyze as a
dynamical system in the plane. Most previous work on Ricci flow on homogeneous
spaces takes the first parametrization (e.g., \cite{IJ}, \cite{KM},
\cite{CS2}). In contrast, we will take the second parametrization, and
consider $\mathcal{M}$ as a quotient of the space of structure constants. This
method has previously been used by the second author to study Ricci flow on
nilmanifolds \cite{P2} (see also \cite{Guz} and \cite{Lau4}). Let $\phi_{t}$
denote the flow on $\mathcal{M}$ determined by the Ricci flow. Theorem A
describes the topological dynamics of the flow $\phi_{t}.$

\begin{figure}[th]
\begin{minipage}[b]{0.4\linewidth}
\centering
\setlength{\unitlength}{2mm} \begin{picture}(40,60)(0,0)
\thinlines
\thicklines
% vertical lines
\put(10,0){\line(0,1){30}}
\put(30,30){\line(0,1){20}}
% horizontal
\put(10,30){\line(1,0){20}}
% slope 1
\put(10,0){\line(2,3){20}}
\put(10,30){\line(1,1){20}}
\put(10,15){\line(4,3){20}}
% dots
\put(10,0){\circle*{1}}
\put(10,15){\circle*{1}}
\put(10,30){\circle*{1}}
\put(30,30){\circle*{1}}
\put(30,50){\circle*{1}}
%  arrows on vert lines
\put(10,8){\vector(0,1){0.8}}
\put(10,22){\vector(0,-1){0.8}}
\put(30,40){\vector(0,1){0.8}}
%  arrows on horizontal lines
\put(20,30){\vector(1,0){0.8}}
% arrows on slanted lines
\put(22,24){\vector(4,3){0.8}}
\put(20,15){\vector(2,3){0.8}}
\put(20,40){\vector(1,1){0.8}}
% point labels
\put(6,0){$p_1$}
\put(6,15){$p_2$}
\put(6,30){$p_1$}
\put(32,30){$p_3$}
\put(32,50){$p_4$}
% trajectory labels
\put(4,22){$T_{1,2}$}
\put(4,8){$T_{1,2}$}
\put(20,12){$T_{1,3}^\prime$}
\put(14,22){$T_{2,3}$}
\put(20,31){$T_{1,3}$}
\put(16,42){$T_{1,4}$}
\put(31,40){$T_{3,4}$}
% basin labels
\put(22,37){$B_{1,4}$}
\put(16,27){$B_{1,3}$}
\put(12,13){$B_{1,3}^\prime$}
\end{picture}
\caption{Schematic of  $\bar{\mathcal{S}}_m$ with Ricci flow lines.}%
\label{figure1}
\end{minipage}\hspace{0.5cm} \begin{minipage}[b]{0.4\linewidth}
\centering
\includegraphics[scale=1]{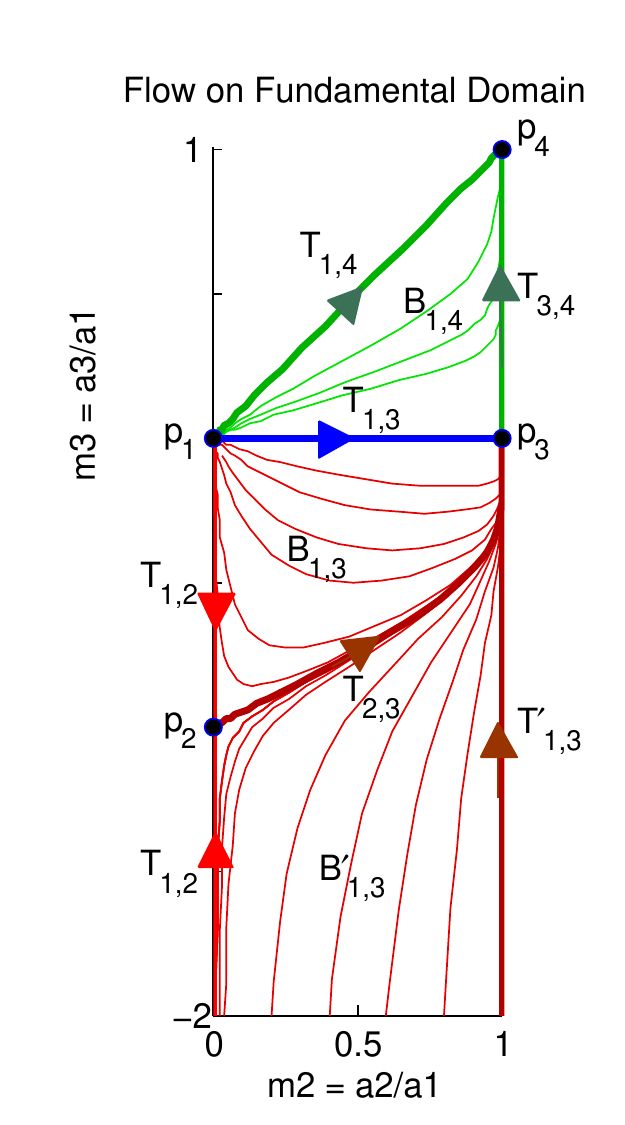}
\caption{The space $\mathcal{S}_m$ with Ricci flow lines.}
\label{figure2}
\end{minipage}\end{figure}

\begin{theorema}
\label{Theorem A}The phase space $\mathcal{M}$ is the disjoint union of the
following invariant sets (see Figures \ref{figure1} and \ref{figure2}):

\begin{itemize}
\item {four points $p_{1},p_{2},p_{3}$ and $p_{4};$}

\item {six one-dimensional trajectories $T_{1,2},T_{1,3},T_{1,3}^{\prime
},T_{1,4},T_{2,3}$ and $T_{3,4}$; and}

\item {three connected two-dimensional open sets $B_{1,4},B_{1,3}$ and
$B_{1,3}^{\prime}${;}}
\end{itemize}

such that

\begin{itemize}
\item {the points $p_{1},p_{2},p_{3}$ and $p_{4}$ are fixed by $\phi_{t};$}

\item {the orbit of a point $p$ in a $T_{i,j}$ or $T_{i,j}^{\prime}$ has
$\lim_{t\rightarrow-\infty}\phi_{t}(p)=p_{i},$ and $\lim_{t\rightarrow\infty
}\phi_{t}(p)=p_{j};$ and }

\item {the orbit of a point $p$ in $B_{i,j}$ or $B_{i,j}^{\prime}$ has
$\lim_{t\rightarrow-\infty}\phi_{t}(p)=p_{i},$ and $\lim_{t\rightarrow\infty
}\phi_{t}(p)=p_{j}.$ }
\end{itemize}
\end{theorema}

Theorem B interprets Theorem A geometrically. In the sequel, we will say a
point in the phase space represents a particular metric, although we actually
mean that it represents the homothety class of the metric, i.e., the
equivalence class of the metric up to isometry and scaling.

\begin{theoremb}
The decomposition of the phase space $\mathcal{M}$ in Theorem A corresponds
geometrically as follows.

\begin{enumerate}
\item {Each of the four fixed points $p_{1},p_{2},p_{3}$ and $p_{4}$
represents a soliton metric: }

\begin{itemize}
\item {$p_{1}$ represents the soliton metric on the three-dimensional
Heisenberg group }$H\left(  3\right)  .$

\item {$p_{2}$ represents the soliton metric on the three-dimensional solvable
group $E(1,1).$}

\item {$p_{3}$ represents the flat metric on the three-dimensional Euclidean
group $\widetilde{E(2)}.$}

\item {$p_{4}$ represents the round metric on the group $\operatorname{SU}%
(2)$.}
\end{itemize}

\item {The five trajectories $T_{1,2},$ $T_{1,3},$ $T_{1,3}^{\prime},$
$T_{1,4},$ and $T_{3,4}$ have Riemannian submersion structures:}

\begin{itemize}
\item $T_{1,2}$ consists of left-invariant metrics on $E\left(  1,1\right)  $
(often denoted $\operatorname{Sol}$). These metrics fiber as Riemannian
submersions over $\mathbb{R}$.

\item $T_{1,3}$ consists of left-invariant metrics on {$\widetilde{E(2)}$}.
These metrics fiber as Riemannian submersions over $\mathbb{R}$.

\item $T_{1,3}^{\prime}$ consists of left-invariant metrics on $\widetilde
{{\operatorname{SL}}_{2}{(\mathbb{R})}}$ which fiber as Riemannian submersions
over the hyperbolic plane $\mathbb{H}^{2}.$

\item $T_{1,4}$ and $T_{3,4}$ consist of left-invariant metrics on
$\operatorname{SU}\left(  2\right)  $ which fiber as Riemannian submersions
over the round sphere $\mathbb{S}^{2}$ (these Riemannian manifolds are often
called Berger spheres). The trajectory $T_{1,4}$ corresponds to submersions
whose fibers are larger than those of the round 3-sphere (corresponding to the
point $p_{4}$) and the trajectory $T_{3,4}$ corresponds to submersions whose
fibers are smaller than those of the round 3-sphere.
\end{itemize}

\item The {three connected open sets $B_{1,4},B_{1,3}$ and $B_{1,3}^{\prime}$
have the structures: }

\begin{itemize}
\item {$B_{1,4}$ consists of left-invariant metrics on $\operatorname{SU}(2)$.
}

\item {$B_{1,3}$ and $B_{1,3}^{\prime}$ consist of left-invariant metrics on
}$\widetilde{{\operatorname{SL}}_{2}{(\mathbb{R})}}${.}
\end{itemize}
\end{enumerate}
\end{theoremb}

Note that the trajectory $T_{2,3}$ is still somewhat mysterious. This
trajectory was discovered independently by Cao, Guckenheimer, and Saloff-Coste
\cite{CGS}, and evidence for it was present in \cite{CNS}. Preliminary
evidence suggests that this trajectory is not invariant under cross curvature
flow, which may indicate it does not arise from extra symmetries of the
Riemannian metric, as the other special orbits do (see Remark
\ref{remark symmetry}).

The organization of this paper is as follows. In Section
\ref{section space of unimod}, we introduce notation and discuss the space of
three-dimensional, unimodular metric Lie algebras and their curvatures. In
Section \ref{section RF} we derive the Ricci flow equations on the space of
structure constants. In Section \ref{dynamics} we analyze the dynamics of the
Ricci flow equations on a natural phase space and then on $\mathcal{M}$,
completing the proof of Theorem A. In Section \ref{geometry section} we
analyze the dynamics geometrically, relating fixed points and special
trajectories to Ricci solitons and Riemannian submersions, proving Theorem B.
In Section \ref{section convergence} we discuss how our convergence results
relate to the current literature on Ricci flows on three-dimensional,
unimodular metric Lie algebras and Lie groups. Finally, we include an appendix
which give the details of characterizing the space $\mathcal{M}$.

\section{Metric Lie algebras and their
curvatures\label{section space of unimod}}

Consider the following definitions.

\begin{definition}
A \emph{metric Lie algebra} $\left(  \mathfrak{g},\mathsf{Q}\right)  $ is a
Lie algebra $\mathfrak{g}$ together with an inner product $\mathsf{Q}$ on
$\mathfrak{g}$. The dimension of the metric Lie algebra is the dimension of
$\mathfrak{g}$, and it is \emph{unimodular} or \emph{nonabelian} if the Lie
algebra $\mathfrak{g}$ is unimodular (i.e., $\operatorname{ad}_{X}$ is trace
free for all $X\in\mathfrak{g}$) or nonabelian, respectively.
\end{definition}

Recall that there is a one-to-one correspondence between Lie algebras
$\mathfrak{g}$ and simply-connected Lie groups $G.$ Furthermore, there is a
one-to-one correspondence between metric Lie algebras $\left(  \mathfrak{g}%
,\mathsf{Q}\right)  $ and simply-connected Lie groups $G$ with a
left-invariant Riemannian metric $g.$ We will use the Riemannian manifold
$\left(  G,g\right)  $ corresponding to the metric Lie algebra $\left(
\mathfrak{g},\mathsf{Q}\right)  $ for the following definition.

\begin{definition}
\label{isometry def}A linear map $L:\mathfrak{g}\rightarrow\mathfrak{g}%
^{\prime}$ is an \emph{isometry} of metric Lie algebras $\left(
\mathfrak{g,}\mathsf{Q}\right)  $ and $\left(  \mathfrak{g}^{\prime
}\mathfrak{,}\mathsf{Q}^{\prime}\right)  $ if it the differential of a
Riemannian isometry $\left(  G,g\right)  \rightarrow\left(  G^{\prime
},g^{\prime}\right)  $ between the corresponding simply-connected Lie groups
with induced left-invariant metrics. Two metrics $\mathsf{Q}$ and
$\mathsf{Q}^{\prime}$ on $\mathfrak{g}$ are \emph{homothetic} if
$\mathsf{Q}^{\prime}=c\mathsf{Q}$ for some $c>0.$ We say $\left(
\mathfrak{g,}\mathsf{Q}\right)  $ and $\left(  \mathfrak{g}^{\prime
}\mathfrak{,}\mathsf{Q}^{\prime}\right)  $ are \emph{equivalent up to isometry
and scaling} if there exists a diffeomorphism $\phi:G\rightarrow G^{\prime}$
and $c>0$ such that $g=c\phi^{\ast}g^{\prime}.$ The space of
three-dimensional, nonabelian, unimodular metric Lie algebras modulo isometry
and scaling will be denoted by $\mathcal{M}$.
\end{definition}

\begin{remark}
\label{remark 1}A related definition is that of isomorphism between metric Lie
algebras. One says two metric Lie algebras $\left(  \mathfrak{g,}%
\mathsf{Q}\right)  $ and $\left(  \mathfrak{g}^{\prime}\mathfrak{,}%
\mathsf{Q}^{\prime}\right)  $ are \emph{isomorphic} if there is a linear map
$L:\mathfrak{g}\rightarrow\mathfrak{g}^{\prime}$ which is a Lie algebra
isomorphism and satisfies $L^{\ast}\mathsf{Q}^{\prime}=\mathsf{Q}.$ The notion
of isomorphism in this definition is stronger than the notion of isometry in
Definition \ref{isometry def}. That is, given an isomorphism of metric Lie
algebras, one can use the group action to extend this to an isometry of the
corresponding simply-connected Lie groups with corresponding left-invariant
metrics. However, it is possible to have isometries of the groups which are
not isomorphisms of metric Lie algebras. For instance, there is a flat metric
$g_{0}$ on the Lie group $E\left(  2\right)  $ of Euclidean transformations
(as well as its universal cover), and there is a flat metric $g_{1}$ on the
abelian group $\mathbb{R}^{3}.$ There is an isometry between $\left(
\widetilde{E\left(  2\right)  },g_{0}\right)  $ and $\left(  \mathbb{R}%
^{3},g_{1}\right)  ,$ but this isometry is not a group isomorphism, and so its
differential is not an isomorphism of metric Lie algebras.
\end{remark}

In the beautiful paper of Milnor \cite{Mil}, geometric properties of
left-invariant metrics on three-dimensional Lie groups are studied in detail.
Milnor computes the curvatures of three-dimensional, unimodular metric Lie algebras:

\begin{theorem}
[{\cite[Theorem 4.3]{Mil}}]\label{mil curvature thm}Suppose $\left(
\mathfrak{g},\mathsf{Q}\right)  $ is a three-dimensional, unimodular metric
Lie algebra. Then there exists a $\mathsf{Q}$-orthonormal frame $\mathcal{B}%
=\{e_{1},e_{2},e_{3}\}$ for $\mathfrak{g}$ such that Lie brackets for
$\mathfrak{g}$ are determined by%
\begin{equation}
\lbrack e_{2},e_{3}]=a_{1}e_{1},\qquad\lbrack e_{3},e_{1}]=a_{2}e_{2}%
,\qquad\text{and}\qquad\lbrack e_{1},e_{2}]=a_{3}e_{3}, \label{diagonalize}%
\end{equation}
for some constants $a_{1},a_{2},a_{3}\in\mathbb{R}.$ Furthermore, this basis
diagonalizes the Ricci endomorphism $\operatorname{Rc}$ such that
\[
\lbrack\operatorname{Rc}]_{\mathcal{B}}=2%
\begin{bmatrix}
\mu_{2}\mu_{3} &  & \\
& \mu_{1}\mu_{3} & \\
&  & \mu_{1}\mu_{2}%
\end{bmatrix}
,
\]
where
\begin{equation}
\mu_{i}=\frac{1}{2}(a_{1}+a_{2}+a_{3})-a_{i}. \label{mu def}%
\end{equation}
The sectional curvatures are given by
\begin{align*}
K(e_{2}\wedge e_{3})  &  =-\mu_{2}\mu_{3}+\mu_{1}\mu_{3}+\mu_{1}\mu_{2}\\
K(e_{3}\wedge e_{1})  &  =\mu_{2}\mu_{3}-\mu_{1}\mu_{3}+\mu_{1}\mu_{2}\\
K(e_{1}\wedge e_{2})  &  =\mu_{2}\mu_{3}+\mu_{1}\mu_{3}-\mu_{1}\mu_{2}.
\end{align*}
Scalar curvature $\rho$ is
\[
\rho=2(\mu_{2}\mu_{3}+\mu_{1}\mu_{3}+\mu_{1}\mu_{2}).
\]

\end{theorem}

Milnor also describes the isomorphism type of a Lie algebra determined by
Equations (\ref{diagonalize}) based on the signs of $a_{1},a_{2},a_{3}$; see
Table \ref{milnor table} (where if the signs are all multiplied by $-1,$ the
Lie algebra is the same).
%TCIMACRO{\TeXButton{B}{\begin{table}[tbp] \centering}}%
%BeginExpansion
\begin{table}[tbp] \centering
%EndExpansion%
\begin{tabular}
[c]{|l|l|l|}\hline
Signs of $\left\{  a_{1},a_{2},a_{3}\right\}  $ & Associated Lie algebra &
Associated Lie groups\\\hline\hline
$+,+,+$ & $\mathfrak{su}\left(  2\right)  \cong\mathfrak{so}\left(  3\right)
$ & $\operatorname{SU}\left(  2\right)  $ or $\operatorname{SO}\left(
3\right)  =\operatorname{Isom}_{+}\left(  \mathbb{S}^{2}\right)  $\\
$+,+,-$ & $\mathfrak{sl}\left(  2,\mathbb{R}\right)  $ & $\operatorname{SL}%
_{2}\left(  \mathbb{R}\right)  $ or $\operatorname{O}\left(  1,2\right)  $ or
$\operatorname{Isom}_{+}\left(  \mathbb{H}^{2}\right)  $\\
$+,+,0$ & $\mathfrak{e}\left(  2\right)  $ & $E\left(  2\right)
=\operatorname{Isom}\left(  \mathbb{E}^{2}\right)  $\\
$+,-,0$ & $\mathfrak{e}\left(  1,1\right)  $ & $E\left(  1,1\right)
=\operatorname{Sol}$\\
$+,0,0$ & $\mathfrak{h}\left(  3\right)  $ & $H\left(  3\right)
=\operatorname{Nil}$\\\hline
\end{tabular}
\caption{The three-dimensional, nonabelian, unimodular Lie groups/algebras.}\label{milnor table}%
%TCIMACRO{\TeXButton{E}{\end{table}}}%
%BeginExpansion
\end{table}%
%EndExpansion

Recalling the space $\mathcal{M}$ of metric Lie algebras from Definition
\ref{isometry def}, according to Theorem \ref{mil curvature thm} we have the
following lemma.

\begin{lemma}
\label{lemma M tilde}Let $\tilde{\Psi}:\left(  \mathbb{R}^{3}\setminus\left\{
\left(  0,0,0\right)  \right\}  \right)  \rightarrow\mathcal{M}$ be the map
which takes $\left(  a_{1},a_{2},a_{3}\right)  $ to the equivalence class of
the metric Lie algebra defined by an orthonormal basis $\left\{  e_{1}%
,e_{2},e_{3}\right\}  $ with Lie bracket determined by (\ref{diagonalize}).
The map $\tilde{\Psi}$ is surjective.
\end{lemma}

Notice that we have excluded the point $(0,0,0),$ representing the abelian Lie
algebra $\mathbb{R}^{3}$, from the domain of $\tilde{\Psi}$. This will allow
$\tilde{\Psi}$ to descend to a map from $\mathbb{RP}^{2}$ to $\mathcal{M}$ so
that we can consider metric Lie algebra equivalence up to scaling.

\begin{proof}
Certainly the above definition defines a Lie bracket, and Theorem
\ref{mil curvature thm} shows that every three-dimensional, unimodular Lie
algebra can be written in this way. Since $\mathcal{M}$ is a quotient of the
set of all such Lie algebras, the result follows.
\end{proof}

We would like to describe the space $\mathcal{M}$ using a fundamental domain.
Since $\tilde{\Psi}$ descends to a map from $\mathbb{RP}^{2}$ to
$\mathcal{M},$ we start with the coordinates $\left(  m_{2},m_{3}\right)
=\left(  \frac{a_{2}}{a_{1}},\frac{a_{3}}{a_{1}}\right)  $ on $\mathbb{RP}%
^{2}$. Let
\begin{equation}
\mathcal{S}_{m}=\left\{  \left(  m_{2},m_{3}\right)  \in\mathbb{R}^{2}:0\leq
m_{2}\leq1\text{ and }m_{3}\leq m_{2}\right\}  , \label{Sm}%
\end{equation}
and let $\sim$ be the equivalence relation on $\mathcal{S}_{m}$ that is
determined by
\begin{equation}
\left(  0,m_{3}\right)  \sim\left(  0,1/m_{3}\right)  \label{tilde def}%
\end{equation}
if $m_{3}\neq0.$ Since a fundamental domain should be compact, we need to
compactify $\mathcal{S}_{m},$ and so we also introduce the compact set
\begin{equation}
\mathcal{\bar{S}}_{m}=\mathcal{S}_{m}\cup\left\{  \infty\right\}
\label{Smbar}%
\end{equation}
which is given the one-point-compactification topology, i.e., open
neighborhoods of $\infty$ consist of the complements of compact subsets of
$\mathcal{S}_{m}.$ We can extend $\sim$ to an equivalence relation on
$\mathcal{\bar{S}}_{m}$ by adding the equivalence
\begin{equation}
\infty\sim\left(  0,0\right)  . \label{tilde def plus}%
\end{equation}
In the appendix, we prove the following:

\begin{theorem}
\label{fundamental domain theorem}There is a bijection between $\mathcal{\bar
{S}}_{m}/\sim$ and $\mathcal{M}$.
\end{theorem}

Using the quotient topology on $\mathcal{\bar{S}}_{m}/\sim,$ there is a
natural topology on $\mathcal{M}$ which makes this map a homeomorphism.

The proof of Theorem \ref{fundamental domain theorem}, as well as some of
discussion in Section \ref{section convergence}, requires the covariant
derivatives of the Ricci tensor, which can be derived in a straightforward way
from the formulas in Theorem \ref{mil curvature thm}. They also appear in
\cite{Las}.

\begin{proposition}
\label{DRc}Suppose $\left(  \mathfrak{g},\mathsf{Q}\right)  $ is a
three-dimensional, unimodular metric Lie algebra. The covariant derivative of
the Ricci operator satisfies%
\[
\left\vert D\operatorname{Rc}\right\vert ^{2}=8\left(  \left(  \mu_{1}-\mu
_{3}\right)  ^{2}\mu_{2}^{4}+\left(  \mu_{1}-\mu_{2}\right)  ^{2}\mu_{3}%
^{4}+\left(  \mu_{2}-\mu_{3}\right)  ^{2}\mu_{1}^{4}\right)  ,
\]
where $\mu_{i}$ are as in Theorem \ref{mil curvature thm}.
\end{proposition}

\section{Ricci deformation of 3D unimodular metric Lie
algebras\label{section RF}}

We now derive the equations for Ricci flow on $\mathcal{M}$. Recall that on a
Riemannian manifold $\left(  M,g\right)  ,$ the Ricci flow is the solution to
the equations
\[
\frac{\partial g}{\partial t}=-2\operatorname{Rc}\left(  g\right)  .
\]
For a left-invariant metric on a Lie group $G$, this flow reduces to a flow of
the inner product $\mathsf{Q}$ on the Lie algebra $\mathfrak{g}$ of $G$; that
is, it reduces to a flow of metric Lie algebras $\left(  \mathfrak{g}%
,\mathsf{Q}_{t}\right)  .$ Recall that Theorem \ref{mil curvature thm} implies
that for any inner product $\mathsf{Q}$ on a three-dimensional, unimodular Lie
algebra, we can find a $\mathsf{Q}$-orthonormal basis $\mathcal{B}=\left\{
e_{1},e_{2},e_{3}\right\}  $ for $\mathfrak{g}$ which diagonalizes the Ricci
tensor. We will see two ways to formulate the Ricci flow:

\begin{enumerate}
\item \label{case 1 evol}Fix the basis $\mathcal{\bar{B}}=\left\{  \bar{e}%
_{1},\bar{e}_{2},\bar{e}_{3}\right\}  $ from Theorem \ref{mil curvature thm}
which is orthonormal with respect to $\mathsf{Q}_{0},$ the initial inner
product, and consider the evolution of the metric coefficients $\mathsf{Q}%
_{t}\left(  \bar{e}_{i},\bar{e}_{j}\right)  .$ In this case, the structure
constants $\lambda_{1},\lambda_{2},\lambda_{3}$ with respect to the basis
$\mathcal{\bar{B}}$, i.e.,
\begin{equation}
\lbrack\bar{e}_{2},\bar{e}_{3}]=\lambda_{1}\bar{e}_{1},\qquad\lbrack\bar
{e}_{3},\bar{e}_{1}]=\lambda_{2}\bar{e}_{2},\qquad\text{and}\qquad\lbrack
\bar{e}_{1},\bar{e}_{2}]=\lambda_{3}\bar{e}_{3}, \label{lambda identity}%
\end{equation}
are fixed in time.

\item \label{case 2 evol}Evolve the basis $\mathcal{B}_{t}=\left\{
e_{1}\left(  t\right)  ,e_{2}\left(  t\right)  ,e_{3}\left(  t\right)
\right\}  $ to be orthonormal with respect to $\mathsf{Q}_{t},$ so that the
metric is the identity in this frame, but the structure constants $a_{1}%
,a_{2},a_{3}$ with respect to $\mathcal{B}_{t}$, which are of the form
(\ref{diagonalize}), depend on time.
\end{enumerate}

In general, if $\mathcal{\bar{B}}=\left\{  \bar{e}_{1},\bar{e}_{2},\bar{e}%
_{3}\right\}  $ is a basis which is orthogonal with respect to a metric
$\mathsf{\bar{Q}}$ and satisfies (\ref{lambda identity}), and
\[
q_{i}=\mathsf{\bar{Q}}\left(  \bar{e}_{i},\bar{e}_{i}\right)
\]
for $i=1,2,3,$ then we see that $\mathcal{B}=\left\{  e_{1},e_{2}%
,e_{3}\right\}  $, where $e_{i}=\bar{e}_{i}/\sqrt{q_{i}},$ is orthonormal and
the structure constants $a_{i}$ (as in Theorem \ref{mil curvature thm}) are
related by
\begin{equation}
a_{i}=\sqrt{\frac{q_{i}}{q_{j}q_{k}}}\lambda_{i}. \label{a lambda corresp}%
\end{equation}
This is how one can relate the solution flows in Formulations
\ref{case 1 evol} and \ref{case 2 evol}. As first described in \cite{IJ}, in
Formulation \ref{case 1 evol} the Ricci curvature is diagonal with respect to
the initial metric $\mathsf{Q}_{0},$ and so the following Ricci flow evolution
can be derived for $q_{1},q_{2},q_{3}$ (see formulas from Theorem
\ref{mil curvature thm}):%
\begin{align}
\frac{d}{dt}\log q_{1}  &  =-4\mu_{2}\mu_{3}=-a_{1}^{2}+(a_{2}-a_{3}%
)^{2}\nonumber\\
\frac{d}{dt}\log q_{2}  &  =-4\mu_{1}\mu_{3}=-a_{2}^{2}+(a_{3}-a_{1}%
)^{2}\label{q-RF equation}\\
\frac{d}{dt}\log q_{3}  &  =-4\mu_{1}\mu_{2}=-a_{3}^{2}+(a_{1}-a_{2}%
)^{2},\nonumber
\end{align}
where the $a_{i}$ are explicit functions of the $q_{i}$ (with fixed parameters
$\lambda_{i}$) defined by (\ref{a lambda corresp}). Thus the equations
(\ref{q-RF equation}) form an autonomous system of ODEs in the variables
$q_{1},q_{2},q_{3}.$ Since the fixed basis $\mathcal{\bar{B}}$ is orthogonal
with respect to $\mathsf{Q}_{t}$ (as determined by $q_{1}\left(  t\right)
,q_{2}\left(  t\right)  ,q_{3}\left(  t\right)  $) and (\ref{lambda identity})
continues to be satisfied at each time $t$, we see that the flow
(\ref{q-RF equation}) really is the Ricci flow for all time. The fact that the
flow remains diagonal is a special property of three-dimensional, unimodular
metric Lie algebras, and is not true in general (see, e.g., \cite{IJL}).

Noting that the right sides of the ODEs (\ref{q-RF equation}) only contain the
$a_{i}$, without explicitly containing the $q_{i}$, it is natural to consider
Formulation \ref{case 2 evol}. The evolution for $a_{i}$ are easily derived
using (\ref{q-RF equation}) and (\ref{a lambda corresp}). Due to Theorem
\ref{fundamental domain theorem} (and the preceding discussion from Section
\ref{section space of unimod}), we will also be interested in
\begin{align}
m_{2}\left(  t\right)   &  =a_{2}\left(  t\right)  /a_{1}\left(  t\right)
,\label{m definition}\\
m_{3}\left(  t\right)   &  =a_{3}\left(  t\right)  /a_{1}\left(  t\right)
.\nonumber
\end{align}
We have the following.

\begin{prop}
\label{flow-equations-unimodular} Let $G$ be a simply-connected,
three-dimensional, nonabelian, unimodular Lie group with left-invariant metric
$g.$ Then the Ricci flow on $G$ with initial metric $g$ corresponds to a flow
of metric Lie algebras $\left(  \mathfrak{g},\mathsf{Q}_{t}\right)  ,$ where
$\mathfrak{g}$ is the Lie algebra of $G$ and $\mathsf{Q}_{0}$ is $g$
restricted to $T_{e}G\cong\mathfrak{g}$. This flow can be realized as a flow
of structure constants $a_{1},a_{2},a_{3}$ (as determined by Theorem
\ref{mil curvature thm}), and, if we suppose that $a_{1}\neq0,$ the ratios
$m_{2}$ and $m_{3}$ (as defined in (\ref{m definition})) obey the equations%
\begin{align}
\frac{dm_{2}}{dt}  &  =a_{1}^{2}m_{2}(1-m_{2})(1+m_{2}-m_{3})
\label{m ode bad}\\
\frac{dm_{3}}{dt}  &  =a_{1}^{2}m_{3}(1-m_{3})(1-m_{2}+m_{3}).\nonumber
\end{align}

\end{prop}

\begin{remark}
We also note that if $a_{1}\neq0,$
\[
\frac{d}{dt}(\log a_{1})=\frac{1}{2}\left[  \frac{d}{dt}(\log q_{1})-\frac
{d}{dt}(\log q_{2})-\frac{d}{dt}(\log q_{3})\right]  =2K(e_{2}\wedge e_{3}).
\]

\end{remark}

The reader may be troubled that the expressions for $dm_{2}/dt$ and
$dm_{3}/dt$ in the proposition are not solely functions of $m_{2}$ and $m_{3}%
$; however, they are useful, because we will be interested in imagining slope
fields for the flow of $m_{2}$ and $m_{3}$ in the $m_{2}$-$m_{3}$ plane. The
common (positive) term $a_{1}^{2}$ simply affects the speed of motion and not
the direction of motion, so the trajectories of $m_{2},m_{3}$ under the ODEs
(\ref{m ode bad}) are the same as the trajectories of the autonomous ODEs
\begin{align}
\frac{dm_{2}}{dt}  &  =m_{2}(1-m_{2})(1+m_{2}-m_{3})\label{RF equation}\\
\frac{dm_{3}}{dt}  &  =m_{3}(1-m_{3})(1-m_{2}+m_{3}).\nonumber
\end{align}

\begin{proof}
The translation to the metric Lie algebra was described at the beginning of
this section. Observe that
\[
m_{2}=\frac{a_{2}}{a_{1}}=\frac{\lambda_{2}}{\lambda_{1}}\sqrt{\frac{q_{2}%
}{q_{1}}}.
\]
Now calculate $dm_{2}/dt$ using (\ref{q-RF equation}) as:
\begin{align*}
\frac{d}{dt}\left(  \frac{a_{2}}{a_{1}}\right)   &  =\frac{\lambda_{2}%
}{\lambda_{1}}\frac{d}{dt}\left(  \sqrt{\frac{q_{2}}{q_{1}}}\right) \\
&  =\frac{\lambda_{2}}{\lambda_{1}}\sqrt{\frac{q_{2}}{q_{1}}}\frac{1}{2}%
\frac{d}{dt}\left(  \log q_{2}-\log q_{1}\right) \\
&  =\frac{1}{2}m_{2}\left(  (a_{3}-a_{1})^{2}-a_{2}^{2}\right)  -\left(
(a_{3}-a_{2})^{2}-a_{1}^{2}\right) \\
&  =m_{2}(a_{1}-a_{2})(a_{1}+a_{2}-a_{3})\\
&  =a_{1}^{2}m_{2}(1-m_{2})(1+m_{2}-m_{3}).
\end{align*}
The formula for $dm_{3}/dt$ follows analogously.
\end{proof}

\begin{remark}
It is clear that the Ricci flow equations (\ref{q-RF equation}) for the
$q_{i}$ determine the $a_{i}$ by (\ref{a lambda corresp}), however one might
ask if the Ricci flow equations for the $a_{i}$ determine the $q_{i}.$ This
is, in fact, true, since once the $a_{i}$ are an explicit function of $t,$ one
can determine the $q_{i}$ by explicitly integrating (\ref{q-RF equation}),
which are now explicit functions of $t.$ This was first observed in \cite{P2}.
\end{remark}

\section{Dynamics of the ODEs\label{dynamics}}

\subsection{Dynamics in $\mathbb{R}^{2}$\label{section R^3 dynamics}}

\begin{figure}[th]
\includegraphics[scale=.8]{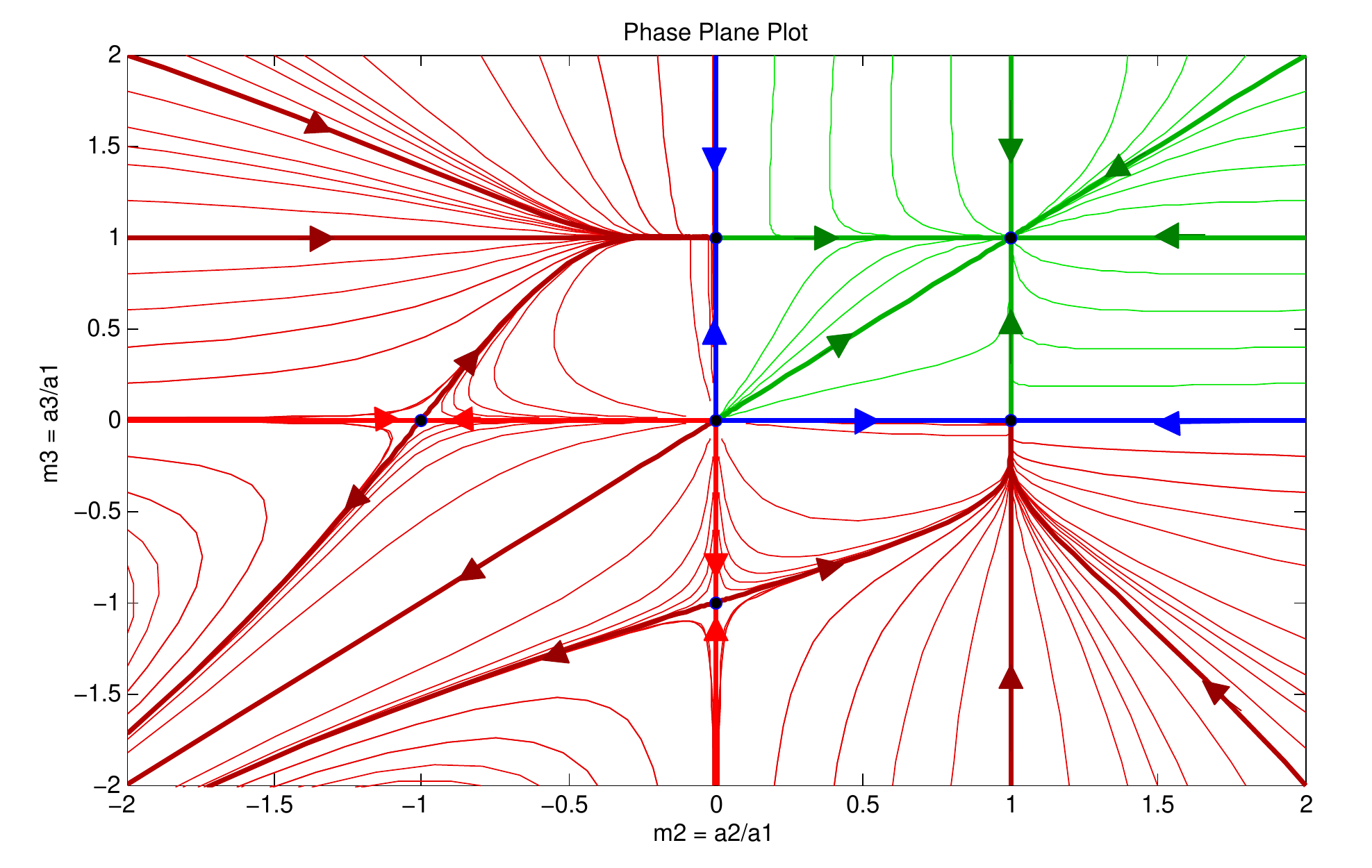}\caption{Phase plane for
ODEs in Equations (\ref{RF equation}).}%
\label{figure3}%
\end{figure}In this section we will look at the qualitative behavior of the
dynamical system (\ref{RF equation}). The phase plane for the system of ODEs
(\ref{RF equation}) is displayed in Figure \ref{figure3}, as computed in
Matlab. First, we consider the fixed points. Since the right sides of each
equation factor into linear terms, it is easy to see that the fixed points of
$\left(  m_{2},m_{3}\right)  $ are $\left(  0,0\right)  ,\left(
0,\pm1\right)  ,\left(  \pm1,0\right)  ,\left(  1,1\right)  .$ Also, it is not
hard to see that the following curves are preserved by the flow: (i)
$m_{2}=0,$ (ii) $m_{3}=0,$ (iii) $m_{2}=1,$ (iv) $m_{3}=1,$ and (v)
$m_{2}=m_{3}.$

The linearization of the system (\ref{RF equation}) is the matrix
\[
\left(
\begin{array}
[c]{cc}%
2m_{2}m_{3}-m_{3}-3m_{2}^{2}+1 & m_{2}\left(  m_{2}-1\right) \\
m_{3}\left(  m_{3}-1\right)  & 2m_{2}m_{3}-m_{2}-3m_{3}^{2}+1
\end{array}
\right)
\]
and so we see that:

\begin{enumerate}
\item $\left(  0,0\right)  $ is an unstable fixed point.

\item $\left(  1,1\right)  $ is a stable fixed point.

\item $\left(  0,-1\right)  $ and $\left(  -1,0\right)  $ are saddle points.
They have stable manifolds tangent to the lines determined by the eigenvectors
$\left(  0,1\right)  $ and $\left(  1,0\right)  $ respectively. They have
unstable manifolds tangent to the lines determined by the eigenvectors
$\left(  2,1\right)  $ and $\left(  1,2\right)  .$

\item $\left(  0,1\right)  $ and $\left(  1,0\right)  $ are degenerate fixed
points. They have stable manifolds tangent to the lines determined by the
eigenvectors $\left(  0,1\right)  $ and $\left(  1,0\right)  $ respectively.
They have a zero eigenvalue corresponding to eigenvectors $\left(  1,0\right)
$ and $\left(  0,1\right)  $; furthermore, one can look in the zero directions
by considering the Taylor series of solutions near $s=0$ for the flow along
curves $s\rightarrow\left(  s,1\right)  $ and $s\rightarrow\left(  1,s\right)
$. We see that, for instance, points near $\left(  1,0\right)  $ and below the
$x$-axis approach the fixed point, while points above the $x$-axis move away
from the fixed point, as seen in Figure \ref{figure3}.
\end{enumerate}

In addition, since $\left(  0,-1\right)  $ and $\left(  -1,0\right)  $ are
saddle points (i.e., the linearizations at these points each have two distinct
eigenvalues of opposite sign), the Stable Manifold Theorem (see, e.g.,
\cite[section 2.7]{Per}) implies each has a one-dimensional unstable manifold.
Although we are unable to calculate the trajectories explicitly, we can, for
instance, compute a Taylor approximation of the curve at $\left(  0,-1\right)
$ to be%
\begin{equation}
m_{3}=-1+\frac{1}{2}m_{2}+\frac{3}{64}m_{2}^{3}+\frac{3}{128}m_{2}^{4}%
+\frac{9}{512}m_{2}^{5}+\frac{57}{4096}m_{2}^{6}+\frac{1461}{131\,072}%
m_{2}^{7}+O\left(  m_{2}^{8}\right)  . \label{taylor series}%
\end{equation}
Furthermore, by the Hartman-Grobman Theorem (\cite[section 2.8]{Per}), the
trajectories of the differential equation in a neighborhood of $\left(
0,-1\right)  $ are homeomorphic to the trajectories of the linearization
around $\left(  0,-1\right)  ,$ and so this curve contains the only trajectory
in the fourth quadrant with $m_{2}\,<1$ which contains $\left(  0,-1\right)
.$

\begin{remark}
We expect that the Taylor series (\ref{taylor series}) has radius of
convergence 1, since the curve has a vertical tangent at the point $\left(
1,0\right)  .$
\end{remark}

\subsection{Dynamics on $\mathcal{M}$}

In this section we prove Theorem A.

\begin{proof}
[Proof of Theorem A]By Theorem \ref{fundamental domain theorem}, we can
restrict our attention to $\mathcal{S}_{m}$ and consider it up to the
equivalence $\sim$ determined by Equation (\ref{tilde def}). The discussion in
Section \ref{section R^3 dynamics} implies that we have the following fixed
points in $\mathcal{S}_{m},$ none of which are equivalent in $\mathcal{\bar
{S}}_{m}$: $p_{1}=\left(  0,0\right)  $, $p_{2}=\left(  0,-1\right)  ,$
$p_{3}=\left(  1,0\right)  ,$ and $p_{4}=\left(  1,1\right)  .$ It is also
clear that for any sequence $\left\{  \left(  m_{2}^{\left(  i\right)  }%
,m_{3}^{\left(  i\right)  }\right)  \right\}  _{i=1}^{\infty}$ with $0\leq
m_{2}^{\left(  i\right)  }\leq1$ and $\lim_{i\rightarrow\infty}m_{3}^{\left(
i\right)  }=-\infty,$ we have
\[
\lim_{i\rightarrow\infty}\left(  m_{2}^{\left(  i\right)  },m_{3}^{\left(
i\right)  }\right)  =\infty\sim p_{1}%
\]
in $\mathcal{\bar{S}}_{m}.$ Define the sets
\begin{align*}
T_{1,2}  &  =\left\{  \left(  m_{2},m_{3}\right)  :m_{2}=0,\ -1<m_{3}%
<0\right\} \\
&  \sim\left\{  \left(  m_{2},m_{3}\right)  :m_{2}=0,\ -\infty<m_{3}%
<-1\right\} \\
{T_{1,3}}  &  {=\{\left(  m_{2},m_{3}\right)  :m_{3}=0,~0<m_{2}<1\}}\\
{T_{1,3}^{\prime}}  &  {=\{\left(  m_{2},m_{3}\right)  :m_{2}=1,~-\infty
<m_{2}<0\}}\\
{T_{1,4}}  &  {=\left\{  \left(  m_{2},m_{3}\right)  :m_{2}=m_{3}%
,~0<m_{2}<1\right\}  }\\
{T_{3,4}}  &  {=\left\{  \left(  m_{2},m_{3}\right)  :m_{2}=1,0<m_{2}%
<1\right\}  .}%
\end{align*}
Notice that these are, in fact, trajectories of the ODEs (\ref{RF equation}),
and correspond to the invariant sets described in Section
\ref{section R^3 dynamics}.

The special trajectory $T_{2,3}$ is defined as the unstable manifold of the
point $p_{2}$ restricted to the set ${\{\left(  m_{2},m_{3}\right)
:0<m_{2}<1\}.}$ As an unstable manifold, it must be invariant. Consider the
set $\left\{  \left(  m_{2},m_{3}\right)  :0<m_{2}<1\text{ and }%
m_{3}<0\right\}  .$ On this set, it follows from (\ref{RF equation}) that if
$1-m_{2}+m_{3}=0$ then $\frac{d}{dt}\left(  m_{2},m_{3}\right)  $ is in the
positive horizontal direction (i.e., $dm_{2}/dt>0$ and $dm_{3}/dt=0$). Thus
the set
\[
\left\{  \left(  m_{2},m_{3}\right)  :0<m_{2}<1\text{, }m_{3}<0\,\text{, and
}1-m_{2}+m_{3}<0\right\}
\]
is invariant. Since on this set we have%
\[
\frac{dm_{3}}{dm_{2}}=\frac{dm_{3}/dt}{dm_{2}/dt}=\frac{m_{3}(1-m_{3}%
)(1-m_{2}+m_{3})}{m_{2}(1-m_{2})(1+m_{2}-m_{3})}>0,
\]
we see that the trajectory $T_{2,3}$ can be written as $m_{3}=f\left(
m_{2}\right)  ,$ where $f$ is a continuous, increasing function for
$0<m_{2}<1$ such that $f\left(  0\right)  =-1$ and $f\left(  1\right)  =0.$

The remaining sets in the partition of $\mathcal{S}_{m}$ are
\begin{align*}
B_{1,4}  &  =\left\{  \left(  m_{2},m_{3}\right)  :0<m_{2}<1,\ 0<m_{3}%
<m_{2}\right\} \\
B_{1,3}  &  =\left\{  \left(  m_{2},m_{3}\right)  :0<m_{2},~f\left(
m_{2}\right)  <m_{3}<0\right\} \\
B_{1,3}^{\prime}  &  =\left\{  \left(  m_{2},m_{3}\right)  :0<m_{2}%
<1,~m_{3}<f\left(  m_{2}\right)  \right\}  .
\end{align*}
These sets are invariant under the flow of the ODEs since their boundaries are
invariant. It remains to show that the sets have the appropriate forward and
backwards limit properties; a straightforward analysis of the phase diagram
completes the proof.{}
\end{proof}

\section{The geometry of the phase space\label{geometry section}}

%geometry, Hopf fibrations, proof of Theorem B
In this section we compare the results of Theorem A with the known geometry of
three-dimensional, unimodular metric Lie algebras, thereby proving Theorem B.

\begin{proof}
[Proof of Theorem B]We first look at the fixed points $p_{1},p_{2},p_{3}%
,p_{4}.$ Recall that for a three-dimensional, unimodular metric Lie algebra
$\left(  \mathfrak{g},\mathsf{Q}\right)  ,$ we have a basis $\left\{
e_{1},e_{2},e_{3}\right\}  $ as described in Theorem \ref{mil curvature thm},
and recall that $q_{i}=\mathsf{Q}\left(  e_{i},e_{i}\right)  .$

\begin{enumerate}
\item $p_{1}=\left(  0,0\right)  .$ This corresponds to the Lie algebra with
structure constants $a_{1}=1$ and $a_{2}=a_{3}=0.$ By Table \ref{milnor table}
we see that this point corresponds to the Heisenberg Lie algebra
$\mathfrak{h}\left(  3\right)  .$ Up to rescaling, there is only one metric
Lie algebra corresponding to $\mathfrak{h}\left(  3\right)  $, and so we see
that $p_{1}$ must be that point. This point corresponds to the Ricci soliton
on $H\left(  3\right)  $ found by Lauret \cite{Lau1}, Baird-Danielo \cite{BD},
and Lott \cite{Lot1}.

\item $p_{2}=\left(  0,-1\right)  .$ This corresponds to the Lie algebra with
structure constants $a_{1}=1$, $a_{2}=0,$ and $a_{3}=-1.$ By Table
\ref{milnor table} we see that $a_{2}=0$ and $a_{3}<0$ determines that this
point corresponds to the solvable Lie algebra $\mathfrak{e}\left(  1,1\right)
.$ Consulting \cite{G2}, we see that the left-invariant Riemannian metrics on
$E\left(  1,1\right)  $ (referred to as $\operatorname{Sol}$ in the reference)
have the form
\begin{equation}
q_{1}\left(  e^{z}dx+e^{-z}dy\right)  ^{2}+q_{2}dz^{2}+q_{3}\left(
e^{z}dx-e^{-z}dy\right)  ^{2}, \label{sol metrics}%
\end{equation}
and the soliton found by Baird-Danielo \cite{BD} and Lott \cite{Lot1} occurs
when $q_{1}=q_{3}.$ Since at the point $p_{2},$ we have $m_{3}=\sqrt
{q_{3}/q_{1}}\lambda_{3}=-1,$ we see that $q_{1}=q_{3},$ and so $p_{2}$
corresponds to this soliton metric. Note that among metrics (\ref{sol metrics}%
), the soliton has the form
\[
2q_{1}\left(  e^{2z}dx^{2}+e^{-2z}dy^{2}\right)  +q_{2}dz^{2},
\]
which has the additional symmetry of switching $x$ and $y.$ Also note that
switching $x$ and $y$ corresponds with switching $q_{1}$ and $q_{3},$ and thus
gives precisely the isometry that identifies the metric Lie algebras
corresponding to the sets
\[
\left\{  \left(  m_{2},m_{3}\right)  :m_{2}=0,~-1<m_{3}<0\right\}
\]
and
\[
\left\{  \left(  m_{2},m_{3}\right)  :m_{2}=0,~-\infty<m_{3}<-1\right\}
\]
when using the equivalence relation on $\mathcal{S}_{m}$ from (\ref{tilde def}).

\item $p_{3}=\left(  1,0\right)  .$ This corresponds to the Lie algebra with
structure constants $a_{1}=a_{2}=1$ and $a_{3}=0.$ By Table \ref{milnor table}
we see that $a_{3}=0$ and $a_{2}>0$ determines that this point corresponds to
the solvable Lie algebra $\mathfrak{e}\left(  2\right)  .$ Consulting
\cite{G2}, we see that the left-invariant Riemannian metrics on $\widetilde
{E\left(  2\right)  }$ (referred to as $\widetilde{\operatorname{Isom}}\left(
\mathbb{E}^{2}\right)  $ in the reference) have the form
\begin{equation}
q_{1}\left(  \sin\theta dx+\cos\theta dy\right)  ^{2}+q_{2}\left(  \cos\theta
dx-\sin\theta dy\right)  ^{2}+q_{3}d\theta^{2}, \label{e2 metrics}%
\end{equation}
with $q_{1}=q_{2}$ determining the flat metric. We see that $p_{3}$
corresponds to the flat metric. Note that the flat metric is the maximally
symmetric metric of the type (\ref{e2 metrics}).

\item $p_{4}=\left(  1,1\right)  .$ This corresponds to the Lie algebra
$\mathfrak{su}\left(  2\right)  $ with $a_{1}=a_{2}=a_{3}=1,$ implying that
$q_{1}=q_{2}=q_{3}.$ One easily sees from Theorem \ref{mil curvature thm} that
this metric has constant sectional curvature, and thus it corresponds to the
round metric on the 3-sphere. Again, we notice that the round metric is
maximally symmetric among all left-invariant metrics on $\operatorname{SU}%
\left(  2\right)  .$
\end{enumerate}

We now look at the special trajectories. It is clear that $T_{1,2}$ and
$T_{1,3}$ correspond to metrics on $E\left(  1,1\right)  $ and $\widetilde
{E\left(  2\right)  }$ respectively, and the explicit metrics shown in
(\ref{sol metrics}) and (\ref{e2 metrics}) show the Riemannian submersion
structures. In \cite{G2}, the left-invariant metrics on $\widetilde
{\operatorname{SL}_{2}\left(  \mathbb{R}\right)  }$ are given explicitly, and
one sees that on the trajectory $T_{1,3}^{\prime}$ we have $m_{2}=-1,$
indicating that $q_{1}=q_{2}$ and that the metrics have the form
\[
q_{1}\frac{1}{y^{2}}\left(  dx^{2}+dy^{2}\right)  +q_{3}\left(  d\theta
-\frac{1}{y}dx\right)  ^{2},
\]
where $\left(  x,y,\theta\right)  \in\mathbb{R\times\mathbb{R}}_{>0}%
\times\mathbb{R}$. These metrics clearly have the form of Riemannian
submersions over the hyperbolic plane.

Now consider the trajectories $T_{1,4}$ and $T_{3,4},$ which we see correspond
to metrics on $\operatorname{SU}\left(  2\right)  .$ Each element in the basis
$\left\{  e_{1},e_{2},e_{3}\right\}  $ exponentiates to a compact group $K$ of
rotations in $\operatorname{SU}\left(  2\right)  .$ The quotient
$\operatorname{SU}\left(  2\right)  /K$ is diffeomorphic to the sphere
$S^{2},$ and the map $\pi:\operatorname{SU}\left(  2\right)  \rightarrow
\operatorname{SU}\left(  2\right)  /K$ is precisely the Hopf fibration (see,
e.g., \cite{Bes}). Using 9.79 and 9.80 in \cite{Bes}, this can be made into a
Riemannian submersion from a left-invariant metric on $\operatorname{SU}%
\left(  2\right)  $ to a $\operatorname{SU}\left(  2\right)  $-invariant
metric on $S^{2}$ with totally geodesic fibers. The Berger spheres are the
metrics on $\operatorname{SU}\left(  2\right)  $ which make $\pi$ a Riemannian
submersion, where $\operatorname{SU}\left(  2\right)  /K$ is given the round
metric on $S^{2}.$ The remaining basis elements span the horizontal subspace
of the submersion, and thus must have equal length for $\pi$ to be a
Riemannian submersion. Thus the submersions are represented only when
$m_{2}=1$ (so $q_{1}=q_{2}$) or $m_{3}=1$ (so $q_{1}=q_{3}$) or $m_{2}=m_{3}$
(so $q_{2}=q_{3}$). The round sphere is when $q_{1}=q_{2}=q_{3}.$ It is now
easy to see that $T_{1,4}$ corresponds to the fibers being larger than the
fibers for the round sphere: $q_{2}=q_{3}$ and $m_{2}<1,$ so we have
$q_{1}>q_{2}=q_{3}.$ Similarly, on $T_{3,4},$ we have $q_{1}=q_{2}$ and
$m_{3}<1,$ so $q_{3}<q_{1}=q_{2}$ and the fibers are smaller than the fibers
of the round sphere.

The descriptions of the basins $B_{1,3},$ $B_{1,3}^{\prime}$, and $B_{1,4}$
follow immediately from Table \ref{milnor table}.
\end{proof}

\begin{remark}
We could have constructed the Riemannian submersions on $\widetilde
{\operatorname{SL}_{2}\left(  \mathbb{R}\right)  }$ in the same way we
constructed the ones for $\operatorname{SU}\left(  2\right)  $ as follows.
Instead of considering $\widetilde{\operatorname{SL}_{2}\left(  \mathbb{R}%
\right)  },$ consider $\operatorname{PSL}_{2}\left(  \mathbb{R}\right)  ,$ the
orientation preserving isometries of $\mathbb{H}^{2}.$ There is a compact
subgroup $K$ acting on the upper half-plane by the isometries
\[
z\rightarrow\frac{\left(  \cos\theta\right)  z-\sin\theta}{\left(  \sin
\theta\right)  z+\cos\theta}%
\]
for all angles $\theta.$ Note that this is the isotropy group of the point
$i.$ Again, by 9.79 and 9.80 in \cite{Bes}, there is a Riemannian submersion
$\pi:\operatorname{PSL}_{2}\left(  \mathbb{R}\right)  \rightarrow
\operatorname{PSL}_{2}\left(  \mathbb{R}\right)  /K$ with totally geodesic
fibers, where $\operatorname{PSL}_{2}\left(  \mathbb{R}\right)  /K$ is given
the geometry of $\mathbb{H}^{2}.$ This submersion may be lifted to the
universal cover to get a line bundle over $\mathbb{H}^{2}.$
\end{remark}

\begin{remark}
\label{remark symmetry}It is an interesting fact that each of the special
points and trajectories except for $T_{2,3}$ correspond to metrics with
additional symmetry. Generically, the left-invariant metrics have isometry
groups of dimension $3.$ The point $p_{2}$ is maximally symmetric among
left-invariant metrics on $E\left(  1,1\right)  $, containing one extra
symmetry, as described in the proof of Theorem B. The points $p_{3}$ and
$p_{4}$ have $6$-dimensional symmetry groups. The trajectories $T_{1,3}%
^{\prime},$ $T_{1,4},$ and $T_{3,4}$ have the additional symmetries of
reparametrizing the fibers: if the fibration structure is locally trivialized
as $\left(  x,s\right)  \in\mathbb{R}^{2}\times S^{1},$ then the map $\left(
x,s\right)  \rightarrow\left(  x,s+\sigma\right)  $ is an isometry. Note that
reparametrizing the fiber is different than multiplication by the generator of
the fiber in the group. These extra symmetries give $T_{1,3}^{\prime},$
$T_{1,4},$ and $T_{3,4}$ four-dimensional isometry groups.
\end{remark}

\section{Remarks on Convergence\label{section convergence}}

In this section, we compare the convergence results of this paper with other
convergence results for Ricci flow on three-dimensional, unimodular Lie
groups. Let $\left(  \mathfrak{g},\mathsf{Q}\right)  $ be a metric Lie
algebra, and let $\left(  G,g\right)  $ be the corresponding simply-connected
Lie group with left-invariant Riemannian metric $g.$ Any Riemannian manifold
$\left(  M,g\right)  $ determines a metric space $\left(  M,d_{g}\right)  ,$
where $d_{g}$ is the induced Riemannian distance function. Often it will be
relevant to consider quotients $G/\Gamma$ which are manifolds, and in the
sequel, we may use $g$ to denote a metric on $G/\Gamma$ as well as on its
universal cover $G.$ The following are all relevant notions of convergence:

\begin{itemize}
\item If the coefficients of the metrics $\mathsf{Q}_{t}$ on the Lie algebra,
which satisfy the Ricci flow ODEs, converge as $t\rightarrow\infty,$ then the
corresponding metrics $g\left(  t\right)  $ converge in $C^{k}$ or $C^{\infty
}$ as tensors. We call this $C^{k}$ or $C^{\infty}$ convergence.

\item The metric spaces $\left(  G/\Gamma,d_{g\left(  t\right)  }\right)  $
converge uniformly as metric spaces (see \cite{CGS}).

\item The metric spaces $\left(  G/\Gamma,d_{g\left(  t\right)  }\right)  $
converge in the pointed Gromov-Hausdorff topology (see \cite{Gro}).

\item The Riemannian groupoids $\left(  G,g\left(  t\right)  ,\Gamma\right)  $
converge in $C^{k}$ or $C^{\infty}$ as Riemannian groupoids (see \cite{Lot1}
and \cite{G2}).
\end{itemize}

Before we describe the previous work, let's consider the convergence in the
present paper. We have looked at convergence of a system of ODEs for a
normalized Ricci flow equation in the space of metric Lie algebras. In
particular, the convergence is for the structure constants $a_{1},a_{2}%
,a_{3},$ which implies convergence of the connection since the connection is
determined by the Lie brackets; in fact, if $D$ is the Riemannian connection
and $\left\{  e_{1},e_{2},e_{3}\right\}  $ is an orthonormal frame as in
Theorem \ref{mil curvature thm}, then
\[
\mu_{i}=\frac{1}{2}\left(  a_{j}+a_{k}-a_{i}\right)  =\left\langle D_{e_{j}%
}e_{k},e_{i}\right\rangle
\]
for $\left\{  i,j,k\right\}  =\left\{  1,2,3\right\}  $ (see (\ref{mu def})).
Convergence in $C^{0}$ of the connections implies convergence in $C^{1}$ of
the Riemannian metrics (see, e.g., \cite[Chapter 3]{Cetal}). The normalization
is not given explicitly, but must be such that none of the Lie bracket
coefficients (for an orthonormal frame) become infinite and the Lie algebra
does not become abelian. This type of convergence is considered on higher
dimensional nilpotent metric Lie algebras in \cite{P2}.

The earliest works consider either the (forward) Ricci flow equation with no
scaling \cite{KM} or the normalized Ricci flow equation where the rescaling is
based on the scalar curvature \cite{IJ}, and prove $C^{0}$ convergence of the
associated left-invariant metrics on the simply-connected Lie group. The
normalized Ricci flow is helpful for $\operatorname{SU}\left(  2\right)  ,$
since it prevents the sphere from shrinking to a point in finite time, but
otherwise the Ricci flows exist for all time $t$ (for both unnormalized and
normalized). Some of these geometries collapse, in the sense that some of the
metric coefficients $q_{1},q_{2},q_{3}$ go to zero, indicating a compact
quotient will have the volume go to zero, as the sectional curvatures go to
zero at a rate of $1/t$. In \cite{CS2}, $C^{0}$ convergence of the backward
direction (of normalized Ricci flow) is considered, and it is found that the
curvatures may go to infinity and convergence is often in a Gromov-Hausdorff
sense to a sub-Riemannian geometry. Similar work is done for the cross
curvature flow in \cite{CNS} and \cite{CS1}. In each of these cases, compact
quotients are considered, and limits are considered collapsing if a
fundamental domain has injectivity radius going to zero as time goes to
infinity. In many cases, compact quotients collapse with bounded curvature,
i.e., the injectivity radius goes to zero while the sectional curvatures stay bounded.

If one is interested in the simply-connected Lie groups, which are
diffeomorphic to $\mathbb{R}^{3}$ in all cases except $\operatorname{SU}%
\left(  2\right)  $, the collapsing does not occur. Instead, one can show
pointed Gromov-Hausdorff convergence to Euclidean space (see \cite{G2}). This
convergence is sometimes due to the fact that the flow is stretching the
geometry, revealing only that the Riemannian manifold is locally Euclidean. To
counteract this effect, Lott \cite{Lot1} considers convergence when the metric
is rescaled by $1/t$ to try to prevent the sectional curvatures from decaying
to zero. In addition, Lott introduces $C^{k}$ convergence on Riemannian
groupoids (see also \cite{G2}), which essentially means that the universal
cover converges in $C^{k}$ while the actions of the fundamental group change
as the metric evolves. In this setting, one finds that the rescaled solutions
on the universal covers converge to expanding soliton metrics, getting results
similar to those in this paper. It would not be hard to carry out a similar
analysis in the backward time direction with the use of \cite{CS2}. If one is
not concerned with a change of topology, the methods of uniform convergence of
metric spaces will suffice, as in \cite{CGS}.

We note that in \cite{Lot1} and \cite{G2}, only $C^{0}$ convergence of the
metrics is shown explicitly, although $C^{\infty}$ convergence of Riemannian
groupoids is claimed. It is also shown that the sectional curvatures converge,
and so using standard compactness arguments (see, for instance, \cite{H3}
\cite{G1} \cite[Chapter 3]{Cetal}), this implies that the convergence is
$C^{1}.$ Under Ricci flow, Shi's work (\cite{Shi1} and \cite{Shi2}) implies
that if the curvature is bounded at $t=0,$ then all covariant derivatives of
curvature converge at all positive times, uniformly for $t\geq\delta>0.$
However, since we are considering rescaling at time $t=0,$ Shi's estimates do
not directly apply (unless we assume uniform backwards existence), so an
additional argument is needed to show that the metrics converge in $C^{k}$ for
$k>1.$ We would like to have at least convergence in $C^{2}$. The cases in
\cite{Lot1} and \cite{G2} can be shown to converge in $C^{k}$ with additional
work bounding the covariant derivative. We do not pursue this here, but do
present an interesting example below concerning a rescaling that does not
produce $C^{\infty}$ convergence.

There is one major difference between the results in this paper and the
results in \cite{Lot1} and \cite{G2}. Metric Lie algebras of type
$\widetilde{\operatorname{SL}_{2}\left(  \mathbb{R}\right)  }$ converge to
flat metrics in our setting, while they converge to $\mathbb{H}^{2}%
\times\mathbb{R}$ in the setting of \cite{Lot1} and \cite{G2}. We note that
$\mathbb{H}^{2}\times\mathbb{R}$ cannot be realized as a three-dimensional,
unimodular Lie group with a left-invariant metric. We see that, if we rescale
the metric to ensure that the sectional curvatures do not go to zero, the Lie
algebra coefficients cannot stay bounded: using the computations in \cite{G2},
we see that
\begin{align*}
q_{1},q_{3}  &  \sim2t,\\
q_{2}  &  \sim E_{1},
\end{align*}
for some constant $E_{1}>0,$ and two sectional curvatures are asymptotically
like $1/t^{2}$ while one is asymptotically like $-1/t.$ The best we can do by
rescaling is to maintain the negative sectional curvature and let the other
two curvatures go to zero. This suggests scaling the metric by $1/t.$ Under
this rescaling, using (\ref{a lambda corresp}) we can calculate that
\begin{align*}
a_{1},a_{3}  &  \sim\left(  \frac{t}{E_{1}}\right)  ^{1/2},\\
a_{2}  &  \sim\frac{1}{2}\left(  \frac{E_{1}}{t}\right)  ^{1/2}.
\end{align*}
Thus we see that this rescaling would not stay in the space of unimodular
metric Lie algebras, and we must take a different scaling to stay in this
space. A different scaling will cause the sectional curvatures all to go to
zero (for instance no scaling at all), revealing Euclidean space.

Metric Lie algebras of the type $\widetilde{E\left(  2\right)  }$ also
converge to a flat metric, so one might ask if rescaling by the maximum
curvature can create a non-flat limit. We will use the notation of \cite{G2}.
Note that%
\[
a_{2}-a_{1}=\frac{q_{2}-q_{1}}{\sqrt{q_{1}q_{2}q_{3}}},
\]
and recall that in this case, Ricci flow has
\begin{align*}
q_{1},q_{2}  &  \sim E_{1}\\
q_{3}  &  \sim E_{3}\\
q_{1}-q_{2}  &  \sim E_{4}e^{-E_{3}t},
\end{align*}
for constants $E_{1},E_{2},E_{3},E_{4},$ and sectional curvatures are all
proportional to $e^{-E_{3}t}.$ If we rescale, replacing $q_{i}$ with
$e^{-E_{3}t}q_{i}$, we see that
\begin{align*}
a_{2}-a_{1}  &  \sim\frac{E_{4}e^{-2E_{3}t}}{E_{5}e^{-\frac{3}{2}E_{3}t}%
}=E_{6}e^{-\frac{1}{2}E_{3}t},\\
a_{1},a_{2}  &  \sim E_{7}e^{\frac{1}{2}E_{3}t},\\
a_{3}  &  =0.
\end{align*}
Thus for the rescaled solution,
\begin{align*}
\mu_{1}  &  \sim\frac{1}{2}E_{6}e^{-\frac{1}{2}E_{3}t}\\
\mu_{2}  &  \sim-\frac{1}{2}E_{6}e^{-\frac{1}{2}E_{3}t}\\
\mu_{3}  &  \sim E_{7}e^{\frac{1}{2}E_{3}t}%
\end{align*}
and thus, using Proposition \ref{DRc},%
\begin{align*}
\left\vert D\operatorname{Rc}\right\vert ^{2}  &  \sim8\left(  \left(
E_{7}e^{\frac{1}{2}E_{3}t}\right)  ^{2}\left(  \frac{1}{2}E_{6}e^{-\frac{1}%
{2}E_{3}t}\right)  ^{4}+\left(  E_{6}e^{-\frac{1}{2}E_{3}t}\right)
^{2}\left(  E_{7}e^{\frac{1}{2}E_{3}t}\right)  ^{4}+\left(  E_{7}e^{\frac
{1}{2}E_{3}t}\right)  ^{2}\left(  \frac{1}{2}E_{6}e^{-\frac{1}{2}E_{3}%
t}\right)  ^{4}\right) \\
&  \sim8E_{6}^{2}E_{7}^{4}e^{E_{3}t},
\end{align*}
and so under this rescaling, $\left\vert D\operatorname{Rc}\right\vert
^{2}\rightarrow\infty$ as $t\rightarrow\infty.$ Since $\left\vert
D\operatorname{Rc}\right\vert ^{2}$ is not bounded, we cannot use
Arzela-Ascoli-type compactness theorems to get convergence of the curvatures
(see \cite{G1} or \cite{Cetal}) and the curvatures may not converge. Thus this
kind of rescaling does not necessarily result in $C^{2}$ and definitely not
$C^{3}$ convergence. Thus there is no natural rescaling which results in a
non-flat limit.

We can also investigate what happens in the backwards time limit and compare
to the results of \cite{CS2}. The results are slightly different, since once
again the rescaling in \cite{CS2} is chosen in a particular way, while we have
chosen a different rescaling. A comparison can be made in a straightforward
way, which we leave to the reader. However, we would like to point out the
particularly interesting case of the backwards limit of the trajectory
$T_{3,4},$ which appears to converge, in our setting, to the point $p_{3}$
(see Figure \ref{figure3}), which represents the flat metric on $\widetilde
{E\left(  2\right)  }$ (see Theorem B). As in the case of convergence of the
forward evolution of $\widetilde{\operatorname{SL}_{2}\left(  \mathbb{R}%
\right)  }$ metrics (i.e., those represented by $B_{1,3},$ $B_{1,3}^{\prime},$
$T_{2,3},$ and $T_{1,3}^{\prime}$ as in Theorem A), the fact that we see
convergence to the flat metric indicates that our rescaling may be revealing
only the infinitesimal Euclidean character of the limit. In fact, as we saw in
Section \ref{geometry section}, the trajectory $T_{3,4}$ corresponds to Berger
metrics which are Riemannian submersions over the round $2$-sphere. As we
follow the trajectory backwards, the fibers shrink and so we expect
convergence to the $2$-sphere. The formalism in \cite{Lot1} and \cite{G2}
together with the calculations in \cite{CS2} allow one to make this
convergence precise in the sense of Riemannian groupoids (and, in particular,
pointed Gromov-Hausdorff).

Finally, we note that the stability analysis of Ricci solitons differs from
that in \cite{IGK}. In \cite{IGK}, only compactly-supported variations are
considered, and it was found, for instance, that the soliton metric on
$H\left(  3\right)  $ is linearly stable. We found that in the space of
three-dimensional, nonabelian, unimodular metric Lie algebras up to scaling,
the soliton metric on $H\left(  3\right)  $ is unstable (it corresponds to a
repulsive fixed point). These two results are not contradictory, since the
variations we consider are certainly not compactly supported.

\appendix{}

\section{Three-dimensional, unimodular metric Lie
algebras\label{appendix unimod}}

In this appendix, we prove a characterization of the space of
three-dimensional, nonabelian, unimodular metric Lie algebras, considered up
to isometry and scaling. Under the correspondence of Lemma \ref{lemma M tilde}%
, three-dimensional, unimodular metric Lie algebras correspond to vectors
$\left(  a_{1},a_{2},a_{3}\right)  \in\mathbb{R}^{3}$ by the formula
(\ref{diagonalize}). To account for equivalence under scaling, we consider
two-dimensional real projective space $\mathbb{RP}^{2}\cong\left(
\mathbb{R}^{3}\setminus\left\{  \left(  0,0,0\right)  \right\}  \right)
{\normalsize /}\mathbb{R}^{\times},$ and denote the image of $\left(
x,y,z\right)  \in\mathbb{R}^{3}\setminus\left\{  \left(  0,0,0\right)
\right\}  $ in $\mathbb{RP}^{2}$ under the quotient map as $\left(
x:y:z\right)  .$ Permuting the components induces an isometry of the metric
Lie algebra, and so we define the action of $\sigma\in S_{3}$, the group of
permutations of three elements, on $\mathbb{RP}^{2}$ by $\sigma\left(
a_{1}:a_{2}:a_{3}\right)  =\left(  a_{\sigma\left(  1\right)  }:a_{\sigma
\left(  2\right)  }:a_{\sigma\left(  3\right)  }\right)  .$ Notice that the
actions of $S_{3}$ and $\mathbb{R}^{\times}$ on $\mathbb{R}^{3}\setminus
\left\{  \left(  0,0,0\right)  \right\}  $ are commutative with respect to
each other, and thus%
\[
\mathbb{RP}^{2}/S_{3}\cong\left(  \mathbb{R}^{3}\setminus\left\{  \left(
0,0,0\right)  \right\}  \right)  {\normalsize /}\left(  \mathbb{R}^{\times
}\times S_{3}\right)  .
\]
We denote the image of $\left(  x,y,z\right)  \in\mathbb{R}^{3}\setminus
\left\{  \left(  0,0,0\right)  \right\}  $ in $\mathbb{RP}^{2}/S_{3}$ as
$\left[  x:y:z\right]  ,$ where the use of square brackets instead of round
brackets indicates the further equivalence using the $S_{3}$ action. We will
now introduce a fundamental domain for $\mathbb{RP}^{2}/S_{3}.$ Recall the
sets $\mathcal{S}_{m}$ and $\mathcal{\bar{S}}_{m}$ defined in (\ref{Sm}) and
(\ref{Smbar}) and the equivalence relation $\sim$ determined by
(\ref{tilde def}) and (\ref{tilde def plus}).

\begin{proposition}
\label{proposition about S}The map%
\[
\tilde{\Phi}:\mathcal{\bar{S}}_{m}\rightarrow\mathbb{RP}^{2}/S_{3}%
\]
defined by%
\begin{align*}
\tilde{\Phi}\left(  m_{2},m_{3}\right)   &  =\left[  1:m_{2}:m_{3}\right]  ,\\
\tilde{\Phi}\left(  \infty\right)   &  =\left[  1:0:0\right]
\end{align*}
is surjective. Moreover, $\tilde{\Phi}$ induces a homeomorphism
\[
\Phi:\mathcal{\bar{S}}_{m}/\sim\rightarrow\mathbb{RP}^{2}/S_{3}.
\]

\end{proposition}

\begin{proof}
Consider a point $\left[  x:y:z\right]  \in\mathbb{RP}^{2}/S_{3}$. We can use
multiplication by $-1$ to ensure that at least one entry is positive and
another is nonnegative. Then we can use the permutations to ensure that $z\leq
y\leq x$, $y\geq0,$ and $x>0.$ In this case we have that
\[
\left[  x:y:z\right]  =\left[  1:\frac{y}{x}:\frac{z}{x}\right]
\]
with
\[
0\leq\frac{y}{x}\leq1
\]
and
\[
\frac{z}{x}\leq\frac{y}{x},
\]
which proves the first statement.

To see that $\Phi$ is well-defined, we must check that $\tilde{\Phi}\left(
0,m_{3}\right)  =\tilde{\Phi}\left(  0,1/m_{3}\right)  $ for $m_{3}\neq0.$ We
see that
\[
\left[  1:0:m_{3}\right]  =\left[  \frac{1}{m_{3}}:0:1\right]  =\left[
1:0:\frac{1}{m_{3}}\right]  .
\]
To see $\Phi$ is a bijection, we note that if
\[
\tilde{\Phi}\left(  m_{2},m_{3}\right)  =\tilde{\Phi}\left(  m_{2}^{\prime
},m_{3}^{\prime}\right)
\]
then
\[
\left[  1:m_{2}:m_{3}\right]  =\left[  1:m_{2}^{\prime}:m_{3}^{\prime}\right]
,
\]
with
\begin{align*}
0  &  \leq m_{2}\leq1,\;\;\;m_{3}\leq m_{2},\\
0  &  \leq m_{2}^{\prime}\leq1,\;\;\;m_{3}^{\prime}\leq m_{2}^{\prime}.
\end{align*}
Thus $\left(  1,m_{2},m_{3}\right)  $ is equal to one of the following, where
$r\neq0$: $r\left(  1,m_{2}^{\prime},m_{3}^{\prime}\right)  $, $r\left(
1,m_{3}^{\prime},m_{2}^{\prime}\right)  $, $r\left(  m_{2}^{\prime}%
,1,m_{3}^{\prime}\right)  $, $r\left(  m_{3}^{\prime},1,m_{2}^{\prime}\right)
$, $r\left(  m_{3}^{\prime},m_{2}^{\prime},1\right)  $, or $r\left(
m_{2}^{\prime},m_{3}^{\prime},1\right)  $. By checking each of the cases, one
can verify that $\Phi$ is a bijection.

To see the continuity, we note that for any sequence $\left(  m_{2}^{\left(
i\right)  },m_{3}^{\left(  i\right)  }\right)  $ satisfying $0\leq
m_{2}^{\left(  i\right)  }\leq1$ and $\lim_{i\rightarrow\infty}m_{3}^{\left(
i\right)  }=-\infty,$ we have
\begin{align*}
\lim_{i\rightarrow\infty}\tilde{\Phi}\left(  m_{2}^{\left(  i\right)  }%
,m_{3}^{\left(  i\right)  }\right)   &  =\lim_{i\rightarrow\infty}\left[
1:m_{2}^{\left(  i\right)  }:m_{3}^{\left(  i\right)  }\right] \\
&  =\lim_{i\rightarrow\infty}\left[  1:\frac{m_{2}^{\left(  i\right)  }}%
{m_{3}^{\left(  i\right)  }}:\frac{1}{m_{3}^{\left(  i\right)  }}\right] \\
&  =\left[  1:0:0\right]  .
\end{align*}

\end{proof}

Certainly there is a map
\[
\Psi:\mathbb{RP}^{2}/S_{3}\rightarrow\mathcal{M}%
\]
induced by the map $\tilde{\Psi}$ defined in Lemma \ref{lemma M tilde}. We
wish to show that this map is bijection. It is immediate that the map is
surjective, so we need only show that it is injective. The main ideas for the
proof are in the paper of Lastaria \cite{Las}, which constructs families of
nonisometric metric Lie algebras with the same Ricci tensors.

We will need two invariants for $\mathcal{M}$. It is easy to see that the
following quantities are invariants of $\mathcal{M}$ if $\left\vert
\operatorname{Rc}\right\vert \neq0$ (and the case $\left\vert
\operatorname{Rc}\right\vert =0$ is distinguished from these cases as well):%
\begin{align*}
\frac{\left(  \lambda_{1}\left[  \operatorname{Rc}\right]  ,\lambda_{2}\left[
\operatorname{Rc}\right]  ,\lambda_{3}\left[  \operatorname{Rc}\right]
\right)  }{\left\vert \operatorname{Rc}\right\vert }  &  \in S^{2}\\
\frac{\left\vert D\operatorname{Rc}\right\vert ^{2}}{\left\vert
\operatorname{Rc}\right\vert ^{3}}  &  \in\mathbb{R}\text{,}%
\end{align*}
where $\lambda_{1}\left[  \operatorname{Rc}\right]  \leq\lambda_{2}\left[
\operatorname{Rc}\right]  \leq\lambda_{3}\left[  \operatorname{Rc}\right]  $
are the eigenvalues of the Ricci operator put into ascending order. Notice
that $\left\vert \operatorname{Rc}\right\vert =0$ only if $\left(  m_{2}%
,m_{3}\right)  =\left(  1,0\right)  ,$ and so the corresponding metric Lie
algebra is not isometric to any metric Lie algebra induced by another element
of $\mathcal{\bar{S}}_{m}/\sim$.

We will consider the map
\[
\tilde{E}=E\circ\Psi\circ\tilde{\Phi}:\mathcal{S}_{m}\setminus\left\{  \left(
1,0\right)  \right\}  \rightarrow S^{2}\times\mathbb{R}%
\]
where
\[
E:\mathcal{M}\rightarrow S^{2}\times\mathbb{R}%
\]
gives the two invariants above. Recall that the Ricci eigenvalues are
\[
\left(  2\mu_{2}\mu_{3},2\mu_{1}\mu_{3},2\mu_{1}\mu_{2}\right)  ,
\]
where $\mu_{1},$ $\mu_{2},$ $\mu_{3}$ are as in Theorem
\ref{mil curvature thm}. We want to express these in terms of $m_{2}$ and
$m_{3},$ so we introduce the following functions which are closely related to
the $\mu_{i}$:
\begin{align}
\nu_{1}  &  =m_{2}+m_{3}-1\nonumber\\
\nu_{2}  &  =1+m_{3}-m_{2}\label{nu def}\\
\nu_{3}  &  =1+m_{2}-m_{3}.\nonumber
\end{align}
Recall that on $\mathcal{S}_{m}$ we have $m_{3}\leq m_{2},$ and so $\nu
_{3}\geq1.$ Also, since $0\leq m_{2}\leq1,$ we have
\begin{equation}
-\nu_{3}\leq\nu_{1}\leq\nu_{2}. \label{nu inequality}%
\end{equation}

We define the partition%
\begin{equation}
P=\left\{  S_{0},S_{1},S_{2},S_{++},S_{-+},S_{++}\right\}
\label{partition def}%
\end{equation}
of $\mathcal{S}_{m}$ by looking at the signs of $\nu_{1}$ and $\nu_{2}$:%
\begin{align*}
S_{0}  &  =\left\{  \left(  1,0\right)  \right\} \\
S_{1}  &  =\left\{  \left(  m_{2},m_{3}\right)  \in\mathcal{S}_{m}:\nu
_{1}=0<\nu_{2}\right\} \\
S_{2}  &  =\left\{  \left(  m_{2},m_{3}\right)  \in\mathcal{S}_{m}:\nu_{1}%
<\nu_{2}=0\right\} \\
S_{++}  &  =\left\{  \left(  m_{2},m_{3}\right)  \in\mathcal{S}_{m}:0<\nu
_{1}\leq\nu_{2}\right\} \\
S_{-+}  &  =\left\{  \left(  m_{2},m_{3}\right)  \in\mathcal{S}_{m}:\nu
_{1}<0<\nu_{2}\right\} \\
S_{--}  &  =\left\{  \left(  m_{2},m_{3}\right)  \in\mathcal{S}_{m}:\nu
_{1}\leq\nu_{2}<0\right\}  .
\end{align*}
Using Theorem \ref{mil curvature thm} and Proposition \ref{DRc} we can write
$\tilde{E},$ which is defined on $\mathcal{S}_{m}\setminus S_{0},$ explicitly
as follows:
\[
\tilde{E}\left(  m_{2},m_{3}\right)  =\left\{
\begin{array}
[c]{cc}%
\left(  \left(  0,0,1\right)  ,\frac{\left(  \nu_{2}^{2}+\nu_{3}^{2}\right)
}{\nu_{2}\nu_{3}}\right)  & \text{if }\left(  m_{2},m_{3}\right)  \in S_{1}\\
\left(  \left(  -1,0,0\right)  ,\frac{\left(  \nu_{1}^{2}+\nu_{3}^{2}\right)
}{\left\vert \nu_{1}\nu_{3}\right\vert }\right)  & \text{if }\left(
m_{2},m_{3}\right)  \in S_{2}\\
\left(  \frac{\left(  \nu_{1}\nu_{2},\nu_{1}\nu_{3},\nu_{2}\nu_{3}\right)
}{\left\vert \nu\nu\right\vert },\ast\right)  & \text{if }\left(  m_{2}%
,m_{3}\right)  \in S_{++}\\
\left(  \frac{\left(  \nu_{1}\nu_{3},\nu_{1}\nu_{2},\nu_{2}\nu_{3}\right)
}{\left\vert \nu\nu\right\vert },\ast\right)  & \text{if }\left(  m_{2}%
,m_{3}\right)  \in S_{-+}\\
\left(  \frac{\nu_{1}\nu_{3},\nu_{2}\nu_{3},\nu_{1}\nu_{2}}{\left\vert \nu
\nu\right\vert },\ast\right)  & \text{if }\left(  m_{2},m_{3}\right)  \in
S_{--},
\end{array}
\right.
\]
where
\[
\ast=\frac{\left(  \nu_{1}-\nu_{3}\right)  ^{2}\nu_{2}^{4}+\left(  \nu_{1}%
-\nu_{2}\right)  ^{2}\nu_{3}^{4}+\left(  \nu_{2}-\nu_{3}\right)  ^{2}\nu
_{1}^{4}}{\left\vert \nu\nu\right\vert ^{3}}%
\]
and
\[
\left\vert \nu\nu\right\vert =\sqrt{\nu_{1}^{2}\nu_{2}^{2}+\nu_{1}^{2}\nu
_{3}^{2}+\nu_{2}^{2}\nu_{3}^{2}}.
\]

Now we will show that $\tilde{E}$ is injective. First we will show that no two
points in different partitions in $P$ correspond to equivalent metric Lie
algebras (equivalent up to isometry and scaling). Then we will show that
within each partition, no two points correspond to equivalent metric Lie algebras.

\begin{proposition}
If $\left(  m_{2},m_{3}\right)  $ and $\left(  m_{2}^{\prime},m_{3}^{\prime
}\right)  $ are in different sets from the partition $P$ (see
(\ref{partition def})), then $\tilde{E}\left(  m_{2},m_{3}\right)  =\tilde
{E}\left(  m_{2}^{\prime},m_{3}^{\prime}\right)  $ implies that

\begin{itemize}
\item $\left(  m_{2},m_{3}\right)  =\left(  0,m_{3}\right)  \in S_{-+}$ and
$\left(  m_{2}^{\prime},m_{3}^{\prime}\right)  =\left(  0,1/m_{3}\right)  \in
S_{--}$ or

\item $\left(  m_{2},m_{3}\right)  =\left(  0,m_{3}\right)  \in S_{--}$ and
$\left(  m_{2}^{\prime},m_{3}^{\prime}\right)  =\left(  0,1/m_{3}\right)  \in
S_{-+}$.
\end{itemize}
\end{proposition}

\begin{proof}
Since $\tilde{E}$ restricted to $S_{++}\cup S_{-+}\cup S_{--}$ does not have
any zeroes in the first component, then certainly we have $\tilde{E}\left(
m_{2},m_{3}\right)  \neq\tilde{E}\left(  m_{2}^{\prime},m_{3}^{\prime}\right)
$ if $\left(  m_{2},m_{3}\right)  $ or $\left(  m_{2}^{\prime},m_{3}^{\prime
}\right)  $ is in $S_{1}$ or $S_{2}.$ Furthermore, since the first component
has all positive entries if $\left(  m_{2},m_{3}\right)  \in S_{++},$ it is
distinguished from $S_{-+}$ and $S_{--},$ both of which have two negative
entries. This shows that $\tilde{E}\left(  m_{2},m_{3}\right)  \neq\tilde
{E}\left(  m_{2}^{\prime},m_{3}^{\prime}\right)  $ for all cases except if
$\left(  m_{2},m_{3}\right)  $ or $\left(  m_{2}^{\prime},m_{3}^{\prime
}\right)  $ in $S_{-+}$ and the other in $S_{--}.$ Say $\left(  m_{2}%
,m_{3}\right)  \in S_{-+}$ and $\left(  m_{2}^{\prime},m_{3}^{\prime}\right)
\in S_{--}$ and $\tilde{E}\left(  m_{2},m_{3}\right)  =\tilde{E}\left(
m_{2}^{\prime},m_{3}^{\prime}\right)  .$ Thus we have that
\[
\frac{\left(  \nu_{1}\nu_{3},\nu_{1}\nu_{2},\nu_{2}\nu_{3}\right)
}{\left\vert \nu\nu\right\vert }=\frac{\left(  \nu_{1}^{\prime}\nu_{3}%
^{\prime},\nu_{2}^{\prime}\nu_{3}^{\prime},\nu_{1}^{\prime}\nu_{2}^{\prime
}\right)  }{\left\vert \nu^{\prime}\nu^{\prime}\right\vert }%
\]
with
\begin{align}
\nu_{1}  &  <0<\nu_{2}\nonumber\\
\nu_{1}^{\prime}  &  \leq\nu_{2}^{\prime}<0\label{nu signs}\\
0  &  <\nu_{3},\nu_{3}^{\prime}.\nonumber
\end{align}
I.e., there exists $c>0$ such that
\begin{equation}
\frac{\nu_{1}\nu_{3}}{\nu_{1}^{\prime}\nu_{3}^{\prime}}=\frac{\nu_{1}\nu_{2}%
}{\nu_{2}^{\prime}\nu_{3}^{\prime}}=\frac{\nu_{2}\nu_{3}}{\nu_{1}^{\prime}%
\nu_{2}^{\prime}}=c, \label{eqn nus ratios}%
\end{equation}
which implies%
\[
\nu_{1}\nu_{2}^{2}\nu_{3}=c^{2}\nu_{1}^{\prime}\left(  \nu_{2}^{\prime
}\right)  ^{2}\nu_{3}^{\prime}.
\]
Using (\ref{eqn nus ratios}), we find that%
\[
\nu_{2}^{2}=c\left(  \nu_{2}^{\prime}\right)  ^{2},
\]
and, using (\ref{nu signs}) and (\ref{eqn nus ratios}) again, we arrive at
\[
\nu_{2}=-\sqrt{c}\nu_{2}^{\prime},\;\;\;\nu_{3}=-\sqrt{c}\nu_{1}^{\prime
},\;\;\;\nu_{1}=-\sqrt{c}\nu_{3}^{\prime}.
\]
It is a consequence of (\ref{nu inequality}) that
\[
-\nu_{3}\leq\nu_{1}=-\sqrt{c}\nu_{3}^{\prime}\leq\sqrt{c}\nu_{1}^{\prime}%
=-\nu_{3},
\]
so
\[
\nu_{1}=-\nu_{3}\text{ and }\nu_{1}^{\prime}=-\nu_{3}^{\prime}.
\]
Thus, using (\ref{nu def}),
\begin{align*}
m_{2}+m_{3}-1  &  =-1-m_{2}+m_{3}\\
m_{2}^{\prime}+m_{3}^{\prime}-1  &  =-1-m_{2}^{\prime}+m_{3}^{\prime},
\end{align*}
or $m_{2}=m_{2}^{\prime}=0.$ Again by (\ref{nu def}), we have the following:
\begin{align*}
\nu_{1}  &  =-\nu_{3}=m_{3}-1\\
\nu_{2}  &  =1+m_{3}\\
\nu_{1}^{\prime}  &  =-\nu_{3}^{\prime}=m_{3}^{\prime}-1\\
\nu_{2}^{\prime}  &  =1+m_{3}^{\prime}.
\end{align*}
Inserting these into (\ref{eqn nus ratios}) gives
\[
\frac{\left(  m_{3}-1\right)  ^{2}}{\left(  m_{3}^{\prime}-1\right)  ^{2}%
}=\frac{m_{3}^{2}-1}{1-\left(  m_{3}^{\prime}\right)  ^{2}},
\]
which implies that
\[
\left(  m_{3}-1\right)  ^{2}\left(  1-\left(  m_{3}^{\prime}\right)
^{2}\right)  -\left(  m_{3}^{2}-1\right)  \left(  m_{3}^{\prime}-1\right)
^{2}=0,
\]
or
\[
\left(  m_{3}-1\right)  \left(  m_{3}^{\prime}-1\right)  \left(  m_{3}%
m_{3}^{\prime}-1\right)  =0.
\]
Since neither $m_{3}$ nor $m_{3}^{\prime}$ can equal one, we must have
$m_{3}m_{3}^{\prime}=1.$
\end{proof}

\begin{proposition}
$\tilde{E}\ $is injective on $S_{1}.$
\end{proposition}

\begin{proof}
On $S_{1},$ we have%
\[
0=\nu_{1}=m_{2}+m_{3}-1,
\]
which implies that
\[
m_{3}=1-m_{2}.
\]
Notice that this implies that
\[
1-m_{2}=m_{3}\leq m_{2},
\]
or%
\[
\frac{1}{2}\leq m_{2}<1
\]
(using that $m_{2}=1$ implies $m_{3}=0).$ Now we can write $\nu_{2}$ and
$\nu_{3}$ in terms of $m_{2}$:
\begin{align*}
\nu_{2}  &  =1+m_{3}-m_{2}=2-2m_{2}\\
\nu_{3}  &  =1+m_{2}-m_{3}=2m_{2}%
\end{align*}
so
\[
\frac{\left(  \nu_{2}^{2}+\nu_{3}^{2}\right)  }{\nu_{2}\nu_{3}}=\frac{\left(
1-m_{2}\right)  ^{2}+m_{2}^{2}}{\left(  1-m_{2}\right)  m_{2}}.
\]
This function is one-to-one on $1/2\leq m_{2}<1.$ Thus $\tilde{E}$ is injective.
\end{proof}

\begin{proposition}
$\tilde{E}\ $is injective on $S_{2}.$
\end{proposition}

\begin{proof}
On $S_{2},$ we have
\[
0=\nu_{2}=1+m_{3}-m_{2},
\]
which implies that
\[
m_{3}=m_{2}-1.
\]
Now we can write $\nu_{1}$ and $\nu_{3}$ in terms of $m_{2}$
\begin{align*}
\nu_{1}  &  =m_{2}+m_{3}-1=2m_{2}-2\\
\nu_{3}  &  =1+m_{2}-m_{3}=2
\end{align*}
so
\[
\frac{\left(  \nu_{1}^{2}+\nu_{3}^{2}\right)  }{\left\vert \nu_{1}\nu
_{3}\right\vert }=\frac{\left(  m_{2}-1\right)  ^{2}+1}{1-m_{2}}.
\]
This function is strictly increasing for $0<m_{2}<1$. The result follows.
\end{proof}

Note that if $U\subset S^{2}\cap\left\{  \left(  x,y,z\right)  :z>0\right\}  $
then the map $\phi:U\rightarrow\mathbb{R}^{2}$ defined by $\phi\left(
x,y,z\right)  =\left(  \frac{x}{z},\frac{y}{z}\right)  $ is one-to-one. We
will use this to prove that $\tilde{E}$ is injective on the sets $S_{++},$
$S_{-+},$ $S_{--}.$

\begin{proposition}
$\tilde{E}\ $is injective on each of the sets $S_{++},$ $S_{-+},$ and
$S_{--}.$
\end{proposition}

\begin{proof}
Let $\pi_{1}$ be the projection onto the first component and let $\bar{E}%
=\phi\circ\pi_{1}\circ\tilde{E}.$ On $S_{++},$ $\bar{E}$ is
\[
\bar{E}\left(  m_{2},m_{3}\right)  =\left(  \frac{\nu_{1}}{\nu_{3}},\frac
{\nu_{1}}{\nu_{2}}\right)  .
\]
We compute the Jacobian determinant:%
\[
\det\left(
\begin{array}
[c]{cc}%
\frac{2\left(  1-m_{3}\right)  }{\nu_{3}^{2}} & \frac{2m_{2}}{\nu_{3}^{2}}\\
\frac{2m_{3}}{\nu_{2}^{2}} & \frac{2\left(  1-m_{2}\right)  }{\nu_{2}^{2}}%
\end{array}
\right)  =\frac{-4\nu_{1}}{\nu_{2}^{2}\nu_{3}^{2}}<0
\]
since $\nu_{1}>0$ on $S_{++}.$ On $S_{-+},$ $\bar{E}$ is
\[
\bar{E}\left(  m_{2},m_{3}\right)  =\left(  \frac{\nu_{1}}{\nu_{2}},\frac
{\nu_{1}}{\nu_{3}}\right)  .
\]
Its Jacobian determinant is
\[
\det\left(
\begin{array}
[c]{cc}%
\frac{2m_{3}}{\nu_{2}^{2}} & \frac{2\left(  1-m_{2}\right)  }{\nu_{2}^{2}}\\
\frac{2\left(  1-m_{3}\right)  }{\nu_{3}^{2}} & \frac{2m_{2}}{\nu_{3}^{2}}%
\end{array}
\right)  =\frac{4\nu_{1}}{\nu_{2}^{2}\nu_{3}^{2}}<0
\]
since $\nu_{1}<0$ on $S_{-+}.$ On $S_{--},$ $\bar{E}$ is
\[
\bar{E}\left(  m_{2},m_{3}\right)  =\left(  \frac{\nu_{3}}{\nu_{2}},\frac
{\nu_{3}}{\nu_{1}}\right)  .
\]
Its Jacobian determinant is
\[
\det\left(
\begin{array}
[c]{cc}%
\frac{2}{\nu_{2}^{2}} & \frac{-2}{\nu_{2}^{2}}\\
\frac{2\left(  m_{3}-1\right)  }{\nu_{1}^{2}} & \frac{-2m_{2}}{\nu_{1}^{2}}%
\end{array}
\right)  =\frac{-4\nu_{3}}{\nu_{1}^{2}\nu_{2}^{2}}>0.
\]
In each case we see that $\bar{E}$ is injective, and thus so is $\tilde{E}.$
\end{proof}

As a consequence of all of these propositions, we can conclude that $\tilde
{E}$ is injective, implying the following theorem, which implies Theorem
\ref{fundamental domain theorem}. Note that we have not yet defined a topology
on $\mathcal{M}$; it will be defined to allow for the theorem.

\begin{theorem}
\label{theorem characterize M}The following spaces are homeomorphic:

\begin{enumerate}
\item $\mathcal{M}$

\item $\mathbb{RP}^{2}/S_{3}$

\item $\mathcal{\bar{S}}_{m}/\sim$
\end{enumerate}
\end{theorem}

\begin{proof}
We have just shown that $\Psi$ gives a bijection between $\mathbb{RP}%
^{2}/S_{3}$ and $\mathcal{M}$, and we can give $\mathcal{M}$ a topology which
makes $\Psi$ a homeomorphism. Proposition \ref{proposition about S} shows that
$\Phi$ is a homeomorphism between $\mathcal{\bar{S}}_{m}/\sim$ and
$\mathbb{RP}^{2}/S_{3}.$
\end{proof}


\begin{thebibliography}{99}                                                                                               %


\bibitem {BD}P. Baird and L. Danielo. Three-dimensional Ricci solitons which
project to surfaces. J. Reine Angew. Math. 608 (2007), 65--91.

\bibitem {BEPS}F. Bourliot, J. Estes, P. M. Petropoulos, and P. Spindel.
Gravitational instantons, self-duality and geometric flows. Preprint at
arXiv:0906.4558v1 [hep-th].

\bibitem {Bes}A. L. Besse. Einstein manifolds. Ergebnisse der Mathematik und
ihrer Grenzgebiete (3) [Results in Mathematics and Related Areas (3)], 10.
Springer-Verlag, Berlin, 1987. xii+510 pp.

\bibitem {CGS}X. Cao, J. Guckenheimer, and L. Saloff-Coste. The backward
behavior of the Ricci and cross curvature flows on SL(2,R). Preprint at
arXiv:0906.4157v1 [math.DG].

\bibitem {CNS}X. Cao, Y. Ni, and L. Saloff-Coste. Cross curvature flow on
locally homogenous three-manifolds. I. Pacific J. Math. 236 (2008), no. 2, 263--281.

\bibitem {CS1}X. Cao and L. Saloff-Coste. Cross curvature flow on locally
homogeneous three-manifolds (II). Preprint at arXiv:0805.3380v1 [math.DG].

\bibitem {CS2}X. Cao and L. Saloff-Coste. Backward Ricci flow on locally
homogeneous three-manifolds. Preprint at arXiv:0810.3352v1 [math.DG].

\bibitem {Cetal}B. Chow, S.-C. Chu, D. Glickenstein, C. Guenther, J. Isenberg,
T. Ivey, D. Knopf, P. Lu, F. Luo, and L. Ni. The Ricci flow: techniques and
applications. Part I. Geometric aspects. Mathematical Surveys and Monographs,
135. American Mathematical Society, Providence, RI, 2007. xxiv+536 pp.

\bibitem {G1}D. Glickenstein. Precompactness of solutions to the Ricci flow in
the absence of injectivity radius estimates. Geom. Topol. 7 (2003), 487--510 (electronic).

\bibitem {G2}D. Glickenstein. Riemannian groupoids and solitons for
three-dimensional homogeneous Ricci and cross-curvature flows. Int. Math. Res.
Not. IMRN 2008, no. 12, Art. ID rnn034, 49 pp.

\bibitem {Gro}M. Gromov. Metric structures for Riemannian and non-Riemannian
spaces. Based on the 1981 French original. With appendices by M. Katz, P.
Pansu and S. Semmes. Translated from the French by Sean Michael Bates. Reprint
of the 2001 English edition. Modern Birkh\"{a}user Classics. Birkh\"{a}user
Boston, Inc., Boston, MA, 2007. xx+585 pp.

\bibitem {IGK}C. Guenther, J. Isenberg, and D. Knopf. Linear stability of
homogeneous Ricci solitons. Int. Math. Res. Not. IMRN 2006, Art. ID 96253, 30 pp.

\bibitem {Guz}G. Guzhvina. Ricci flow on almost flat manifolds. Differential
geometry and its applications, 133--146, World Sci. Publ., Hackensack, NJ, 2008.

\bibitem {H1}R. S. Hamilton. Three-manifolds with positive Ricci curvature. J.
Differential Geom. 17 (1982), no. 2, 255--306.

\bibitem {H2}R. S. Hamilton. The formation of singularities in the Ricci flow.
Surveys in differential geometry, Vol. II (Cambridge, MA, 1993), 7--136, Int.
Press, Cambridge, MA, 1995.

\bibitem {H3}R. S. Hamilton. A compactness property for solutions of the Ricci
flow. Amer. J. Math. 117 (1995), no. 3, 545--572.

\bibitem {IJ}J. Isenberg and M. Jackson. Ricci flow of locally homogeneous
geometries on closed manifolds. J. Differential Geom. 35 (1992), no. 3, 723--741.

\bibitem {IJL}J. Isenberg, M. Jackson, and P. Lu. Ricci flow on locally
homogeneous closed 4-manifolds. Comm. Anal. Geom. 14 (2006), no. 2, 345--386.

\bibitem {Jen}G. R. Jensen. The scalar curvature of left-invariant Riemannian
metrics. Indiana Univ. Math. J. 20 1970/1971 1125--1144.

\bibitem {KM}D. Knopf and K. McLeod. Quasi-convergence of model geometries
under the Ricci flow. Comm. Anal. Geom. 9 (2001), no. 4, 879--919.

\bibitem {Las}F. G. Lastaria. Homogeneous metrics with the same curvature.
Simon Stevin 65 (1991), no. 3-4, 267--281.

\bibitem {Lau1}J. Lauret. Ricci soliton homogeneous nilmanifolds. Math. Ann.
319 (2001), no. 4, 715--733.

\bibitem {Lau2}J. Lauret. Degenerations of Lie algebras and geometry of Lie
groups. Differential Geom. Appl. 18 (2003), no. 2, 177--194.

\bibitem {Lau3}J. Lauret. Einstein solvmanifolds and nilsolitons. Preprint at
arXiv:0806.0035v1 [math.DG].

\bibitem {Lau4}J. Lauret. Homogeneous Ricci flows and solitons. Preprint.

\bibitem {Lot1}J. Lott. On the long-time behavior of type-III Ricci flow
solutions. Math. Ann. 339 (2007), 627-666.

\bibitem {Lot2}J. Lott. Dimensional reduction and the long-time behavior of
Ricci flow. Preprint at arXiv:0711.4063v3 [math.DG].

\bibitem {Mil}J. Milnor. Curvatures of left invariant metrics on Lie groups.
Advances in Math. 21 (1976), no. 3, 293--329.

\bibitem {P1}T. Payne. The existence of soliton metrics for nilpotent Lie
groups. Geom. Dedicata (2009), published online.

\bibitem {P2}T. Payne. The Ricci flow for nilmanifolds. Preprint at
arXiv:0812.2229v1 [math.DG].

\bibitem {Per}L. Perko. Differential equations and dynamical systems. Second
edition. Texts in Applied Mathematics, 7. Springer-Verlag, New York, 1996.
xiv+519 pp.

\bibitem {Shi1}W.-X. Shi. Ricci deformation of the metric on complete
noncompact Riemannian manifolds. J. Differential Geom. 30 (1989), no. 2, 303--394.

\bibitem {Shi2}W.-X. Shi. Complete noncompact three-manifolds with nonnegative
Ricci curvature. J. Differential Geom. 29 (1989), no. 2, 353--360.
\end{thebibliography}
\end{document}